\documentclass[11pt]{amsart}
\usepackage{amssymb, epsf, colordvi}
\usepackage{multirow}
\usepackage{graphicx}
\usepackage[all]{xy}
\newif\ifEPSF
\EPSFtrue %
\numberwithin{equation}{section}

\newtheorem{thm}{Theorem}
\numberwithin{thm}{section}
\newtheorem{prop}[thm]{Proposition}
\newtheorem{lemma}[thm]{Lemma}
\newtheorem{cor}[thm]{Corollary}

\newtheorem{example}[thm]{Example}
\newtheorem{remark}[thm]{Remark}
\newtheorem{definition}[thm]{Definition}

\newenvironment{rem}{\begin{remark}\rm}{\end{remark}}
\newcounter{FNC}[page]
\def\fauxfootnote#1{{\addtocounter{FNC}{1}$^\fnsymbol{FNC}$%
     \let\thefootnote\relax\footnotetext{$^\fnsymbol{FNC}$#1}}}

\copyrightinfo{}{}

\newcommand{\calS}{\mathcal{S}}

\newcommand{\calM}{\mathcal{M}}

\newcommand{\C}{\mathbb{C}}
\newcommand{\T}{\mathbb{T}}
\newcommand{\Trop}{\mathbb{T}{\! {\tiny \mbox{rop}}}}
\newcommand{\Rtrop}{\mathbb{R}{\! {\tiny \mbox{trop}}}}

\newcommand{\R}{\mathbb{R}}

\newcommand{\Z}{\mathbb{Z}}
\newcommand{\rk}{\mbox{rk} \,}

\newcommand\CQFD{\hfill $\Box$ \newline}

\newcommand\KK{{\mathbb{K}}}
\newcommand\ZZ{{\mathbb{Z}}}

\newcommand\RR{{\mathbb{R}}}
\newcommand\CC{{\mathbb{C}}}

\newcommand\bin[2]{{\mathrm{C}^{#2}_{#1}}}

\newcommand\vol{\operatorname{vol}}

\newcommand\cal[1]{\mathcal{#1}}

\newcommand\val{\operatorname{val}}

\newcommand\codim{\operatorname{codim}}

\newcommand\trop{\operatorname{trop}}

\newcommand\VV{W}
\newcommand\VC{V_\mathbb{C}}

\DeclareMathOperator{\Hom}{Hom}
\DeclareMathOperator{\Spec}{Spec}

\title{Euler Characteristic of real nondegenerate tropical complete intersections}

\author{Benoit Bertrand}
\address{Section de math\'ematiques\\
Universit\'e de Gen\`eve\\
case postale 64\\
2-4 rue du li\`evre\\
1211 Gen\`eve 4\\
Suisse}
\email{benoit.bertrand@math.unige.ch}
\urladdr{http://www.unige.ch/math/folks/bertrand}

\author{Frederic Bihan}
\address{Laboratoire de Math\'ematiques\\
         Universit\'e de Savoie\\
         73376 Le Bourget-du-Lac Cedex\\
         France}

\email{Frederic.Bihan@univ-savoie.fr}
\urladdr{http://www.lama.univ-savoie.fr/~bihan}

\thanks{Bertrand was partially supported by the European research
  network IHP-RAAG contract HPRN-CT-2001-00271 and whishes to thank
  Max Planck Institut f\"ur Mathematik for excellent working conditions}

\begin{document}


\begin{abstract}
We define nondegenerate tropical complete intersections
imitating the corresponding definition in complex algebraic geometry.
As in the complex situation, all nonzero intersection multiplicity numbers
between tropical hypersurfaces defining a nondegenerate tropical complete intersection
are equal to $1$. The intersection multiplicity numbers we use are
sums of mixed volumes of polytopes which are dual to cells of the tropical hypersurfaces.
We show that the Euler characteristic of a real nondegenerate tropical complete intersection
depends only on the Newton polytopes of the tropical polynomials which define
the intersection.  Basically, it is equal to the usual signature of
a complex complete intersection with same Newton polytopes, when this signature is defined.
The proof reduces to
the toric hypersurface case, and uses the notion of $E$-polynomials of
complex varieties.
\end{abstract}

\maketitle

%
\section*{Introduction}

Tropical geometry appeared recently in various fields of mathematics
(See \cite{RGST}, \cite{EKL04},  \cite{PaSt04}, \cite{MCounting}).
Tropical varieties can be defined as
the topological closure of the image under the valuation of algebraic
varieties over the field of Puiseux Series $\KK$. For example $T$ is a
tropical hypersurface if there exists an algebraic hypersurface
$Z_\KK$ in $(\KK^*)^n$ such that $T=\overline{V(Z)}$ where $V$ is the
coordinatewise valuation (we rather take minus the valuation).  By a
theorem due to Kapranov (see Theorem~\ref{T:Kapranov}) tropical
hypersurfaces are nonlinearity loci of piecewise-linear convex
functions on $\RR^n$ of the form $f^{\trop}(x)= \max_{\omega \in
\Omega}(<x,\omega>-a_\omega)$ where $\Omega$ is a finite subset of
$\ZZ^n$ and $a_\omega$ a real number.  One of the important
application of tropical geometry is due to Mikhalkin~\cite{MCounting}
who gave a combinatorial way to count the number of curves of given
degree and genus passing through the appropriate number of given
generic points.  Mikhalkin's proof uses a complexification of tropical
curves and a patchworking principle. The algorithm exposed in
\cite{MCounting} has a real counterpart for which it is necessary to
introduce the real part of complexified tropical curves.  The relation
between these real tropical objects and objects appearing in Viro
combinatorial patchworking method is very deep (see ~\cite{Vir01}, for
example).  Actually nonsingular real tropical hypersurfaces are
equivalent from the topological point of view to the so called
primitive $T$-hypersurfaces appearing in the combinatorial Viro
method.

In \cite{Stuint}, Bernd Sturmfels generalized the combinatorial
patchworking method to complete intersections (see
Section~\ref{S:tcon}). The above definition also applies for tropical
varieties.  Namely, one can define a tropical variety to be the image
of an algebraic variety over $\KK$ under the valuation map (see
Section~\ref{S:standard}). This leads also to the notion of complex
tropical variety and real tropical variety (see
Section~\ref{cxrealtropvar}).  In Section~\ref{S:nondeg} we give a
definition for the notion of nondegenerate tropical complete
intersection which builds upon the definition of a nonsingular tropical hypersurface in a manner
similar to the classical complex situation, recalled in
Section~\ref{S:toric}.  We extend in Section~\ref{S:mult} the
definition of tropical intersection multiplicity numbers which was
introduced by Mikhalkin in~\cite{Mtropappli} and show that our
definition is consistent with the classical situation. In particular,
all intersection multiplicity numbers which occur in a nondegenerate
tropical complete intersection are equal to $1$ (or $0$).  We think that our
definition of tropical intersection multiplicity numbers can be of
independent interest.  The goal of this paper is to extend a previous
result of the first author (see~\cite{Be1}) from the case of hypersurfaces to the case of
complete intersections. Roughly speaking, we prove that if
$f_1,\dotsc,f_k$ are polynomials in $\KK[z_1,\dotsc,z_n]$ which define
a nondegenerate tropical complete intersection $Y^{\trop}$, then the Euler
characteristic of the corresponding real tropical complete intersection $\R
Y^{\trop}$ depends only on the Newton polytopes
$\Delta_1,\dotsc,\Delta_k$ of the polynomials and is equal to the
mixed signature $\tilde{\sigma}(Y)$ of a generic complex intersection
$Y$ defined by complex polynomials with the same Newton polytopes.
The precise statement is given in Theorem~\ref{maintorus}.  The notion
of mixed signature is defined by means of the so called
$E$-polynomials (see Section~\ref{S:Epol}).  When $Y$ is a projective
complete intersection of even dimension (over $\C$), then the mixed
signature $\tilde{\sigma}(Y)$ is equal to the usual signature
$\sigma(Y) = \sum_{p+q=0\> [2]} {(-1)}^p h^{p,q}(Y)$, where the
$h^{p,q}(Y)$ are Hodge numbers.  One advantage of the mixed signature
is that it is defined even for a non-projective variety and is additive, as
it is the case for the Euler characteristic.  With the help of this
additivity property, we are able to reduce the proof of the main
result to a proof of the toric hypersurface case. The proof of the
toric hypersurface case uses heavily results obtained by V. Batyrev
and L. Borisov in the paper~\cite{Ba-Bo}.

\section{Toric geometry}\label{S:toric}

We fix some notations and recall some standard properties of toric
geometry.  We refer to~\cite{Ful} for more details.  Let $N \simeq
\Z^n$ be a lattice of rank $n$ and $M=\Hom_{\Z}(N,\Z)$ be its dual
lattice. The associated complex torus is
$\T_{N}:=\Spec(\C[M])=\Hom_{\Z}(M,\C^*)=N \otimes_{\Z} \C^*\simeq
(\C^*)^n$.  Let $f \in {\C} [M ]$ be a Laurent polynomial in the group
algebra associated with $M$
$$f(x)=\sum c_m x^m,$$ where each $m$ belongs to $M$ and only a finite
number of $c_m$ are nonzero.  We will usually have $M={\Z}^n$, so
that ${\C} [M ]=\C[x_1^{\pm 1}, \dotsc, x_n^{\pm 1}]$.  The {\it
support} of $f$ is the subset of $M$ consisting of all $m$ such that
the coefficient $c_m$ is nonzero.  The convex hull of this support in
the real affine space generated by $M$ is called {\it the Newton
polytope} of $f$. This is a {\it lattice} polytope, or a polytope {\it
with integer vertices}, which means that all the vertices of $\Delta$
belong to $M$. In this paper all polytopes will be lattice polytopes
and the ambient lattice $M$ will be clear from the context. We denote
by $M(\Delta)$ the saturated sublattice of $M$ which consists of all
integer vectors parallel to $\Delta$ and by $N(\Delta)$ the dual
lattice. The dimension of $\Delta$ is the rank of $M(\Delta)$, or
equivalently the dimension of the real vector space $M(\Delta)_{\R}$
generated by $\Delta$.  The polynomial $f$ (or rather $x^{-m}f \in
\C[M(\Delta)]$ for any choice of $m$ in the support of $f$) defines an
hypersurface $Z_f$ in the torus $\T_{N(\Delta)}$.  Let $X_{\Delta}$
denote the projective toric variety associated with $\Delta$. The
variety $X_{\Delta}$ contains $\T_{N(\Delta)}$ as a dense Zarisky open
subset and we denote by $\bar{Z}_f$ the Zarisky closure of $Z_f$ in
$X_{\Delta}$.  Let $\Gamma$ be any face of $\Delta$. If $f^{\Gamma}$
is {\it the truncation} of $f$ to $\Gamma$, that is, the polynomial
obtained from $f$ by keeping only those monomials whose exponents
belong to $\Gamma$, then $\bar{Z}_f \cap
\T_{N(\Gamma)}=Z_{f^{\Gamma}}$ and $\bar{Z}_f \cap
X_{\Gamma}=\bar{Z}_{f^{\Gamma}}$.  We have the classical notion of nondegenerate Laurent polynomial.

\begin{definition}\label{D:classicalnondegenerate}
A polynomial $f$ with Newton polytope $\Delta$ is called {\it nondegenerate} if for any face $\Gamma$ of $\Delta$ of positive dimension
(including $\Delta$ itself), the toric hypersurface $Z_{f^{\Gamma}}$
is a nonsingular hypersurface.
\end{definition}
Note that if $\Gamma$ is a vertex of $\Delta$, then $Z_{f^{\Gamma}}$
is empty. In the previous definition, one may equivalently consider
$f^{\Gamma}$ as a polynomial in $\C[M]$ and thus look at the
corresponding hypersurface of the whole torus $\T_{N}$. Indeed, this
hypersurface of $\T_N$ is the product of $Z_{f^{\Gamma}} \subset
\T_{N(\Gamma)}$ with the subtorus of $\T_N$ corresponding to a
complement of $M(\Gamma)$ in $M$.  If $\Delta$ is the Newton polytope
of $f$, then the projective hypersurface $\bar{Z}_f \subset
X_{\Delta}$ is nonsingular if and only if $f$ is nondegenerate and
$X_{\Delta}$ has eventually a finite number of singularities which are
zero-dimensional $\T_{N(\Delta)}$-orbits corresponding to vertices of
$\Delta$.
Consider polynomials $f_1,\dotsc,f_k \in {\C} [M ]$
and denote by $\Delta_i$ the Newton polytope of $f_i$.
Let $\Delta$ be the {\it Minkowsky sum} of these polytopes
$$\Delta=\Delta_1+\cdots+\Delta_k.$$

Each polynomial $f_i$ seen as a polynomial in $\C[M(\Delta)]$ defines
a toric hypersurface $Z_{f_i,\Delta}$ in $\T_{N(\Delta)}$ and it makes
sense to consider the toric intersection
\begin{equation}\label{E:toricintersection}
Z_{f_1,\Delta} \cap \cdots \cap Z_{f_k,\Delta} \subset \T_{N(\Delta)}.
\end{equation}
Denote by $\bar{Z}_{f_i,\Delta}$ the Zarisky closure in $X_{\Delta}$
of $Z_{f_i,\Delta}$.  For each $i=1,\dotsc,k$ there is a toric
surjective map $\rho_i: X_{\Delta} \rightarrow X_{\Delta_i}$ such that
$Z_{f_i,\Delta} = \rho_i^{-1}({Z}_{f_i})$ and $\bar{Z}_{f_i,\Delta} =
\rho_i^{-1}(\bar{Z}_{f_i})$.  This leads to
\begin{equation}\label{E:projectiveintersection}
 \bar{Z}_{f_1,\Delta} \cap \cdots \cap \bar{Z}_{f_k,\Delta} \subset X_{\Delta} .
\end{equation}
Each face $\Gamma$ of $\Delta$ can be uniquely written as a Minkowsky sum
\begin{equation}\label{E:Minkowskyface}
\Gamma=\Gamma_1+\cdots+\Gamma_k
\end{equation}
where $\Gamma_i$ is a face of $\Delta_i$. Substituting the truncation $g_i:=f_i^{\Gamma_i}$ to $f_i$ and $\Gamma_i$
to $\Delta_i$ gives the toric intersection
\begin{equation}\label{E:toricintersectionbis}
Z_{g_1,\Gamma} \cap \cdots \cap Z_{g_k,\Gamma} \subset \T_{N(\Gamma)}.
\end{equation}
which leads to
\begin{equation}\label{E:projectiveintersectionbis}
\bar{Z}_{g_1,\Gamma} \cap \cdots \cap  \bar{Z}_{g_k,\Gamma} \subset X_{\Gamma} .
\end{equation}
Similarly to the hypersurface case the intersection
of~\eqref{E:projectiveintersection} with $\T_{N(\Gamma)}$ (resp.,
$X_{\Gamma}$) coincides with~\eqref{E:toricintersectionbis}
(resp.,~\eqref{E:projectiveintersectionbis}).  Moreover, the
intersection~\eqref{E:projectiveintersection} is the union over all
faces $\Gamma$ of $\Delta$ of the toric
intersections~\eqref{E:toricintersectionbis}.

The {\it Cayley polynomial} associated with $f_1,\dotsc,f_k$ is the
polynomial $F \in \C[M \oplus {\Z}^k]$ defined by
\begin{equation}\label{E:Cayleypolynomial}
F(x,y)=\sum_{i=1}^k y_if_i(x).
\end{equation}
Its Newton polytope is the {\it Cayley polytope} associated with
$\Delta_1,\dotsc,\Delta_k$ and will be denoted by
\begin{equation}\label{E:Cayleypolytope}
C(\Delta_1,\dotsc,\Delta_k) \subset M_{\R} \times {\R}^k.
\end{equation}
Since $F$ is a homogeneous (of degree $1$) with respect to the
variable $y$, the polytope $C(\Delta_1,\dotsc,\Delta_k)$ lies on a
hyperplane and has thus dimension at most $n+k-1$. In fact, the
dimension of $C(\Delta_1,\dotsc,\Delta_k)$ is $\dim(\Delta)+k-1$.  The
faces of $C(\Delta_1,\dotsc,\Delta_k)$ are themselves Cayley
polytopes.  Namely, the faces of $C(\Delta_1,\dotsc,\Delta_k)$ are the
Newton polytopes of all polynomials
$$\sum_{i \in I} y_if_i^{\Gamma_i}(x)$$ such that $\emptyset \neq I
\subset \{1,\dotsc,k\}$ and $\Gamma=\sum_{i \in I} \Gamma_i$ is a face
of $\sum_{i \in I} \Delta_i$ with $\Gamma_i$ a face of $\Delta_i$ for
each $i$.  We will call {\it admissible} such a collection
$(\Gamma_i)_{i \in I}$. Note that by face we do not mean proper face.
In particular $(\Delta_i)_{i \in I}$ is admissible for any non empty
subset $I$ of $\{1,\dotsc,k\}$.  If $(\Gamma_i)_{i \in I}$ is
admissible, we also call admissible the collection of polynomials
$(f_i^{\Gamma_i})_{i \in I}$ and the corresponding toric intersection
\begin{equation}\label{E:toricintersectionbisfaces}
\bigcap_{i \in I} Z_{f_i^{\Gamma_i}, \Gamma} \subset \T_{N(\Gamma)} .
\end{equation}

\begin{definition}\label{D:nondegenerateclassicalinters}
The collection $(f_1,\dotsc,f_k)$ is {\it nondegenerate} if the
associated Cayley polynomial $F(x,y)=\sum_{i=1}^k y_if_i(x)$ is nondegenerate.
\end{definition}

The following result is based on the classical {\it Cayley trick}
(see, for example, \cite{GKZ}).

\begin{prop}\label{P:Cayleytrick}
The collection $(f_1,\dotsc,f_k)$ is nondegenerate if and only if any
admissible toric intersection ~\eqref{E:toricintersectionbisfaces} is
a complete intersection.
\end{prop}
\begin{proof}
As mentionned earlier, we can consider the polynomials
$f_i^{\Gamma_i}$ occurring in~\eqref{E:toricintersectionbisfaces} as
polynomials in $\C[M]$ and thus look at the corresponding intersection
in the whole torus $\T_N$.  An easy computation shows that if hypersurfaces
defined by polynomials $g_i \in \C[M]$, $i \in I$, do not intersect
transversally at a point $X \in \T_N$, then there exists
$\lambda=(\lambda_j )_{j \in J} \in (\C^*)^{|J|}$ with $J \subset I$
so that $\sum_{j \in J} y_j g_j(x)$ defines an hypersurface with a
singular point at $(X,\lambda) \in \T_N \times (\C^*)^{|J|}$.
Similarly, if a truncation $\sum_{i \in I} y_ig_i(x)$ of $F$ to a face
of $C(\Delta_1,\dotsc,\Delta_k)$ defines an hypersurface with a
singular point $(X,\lambda)$ in the corresponding torus, then the
hypersurfaces defined by $g_i$ for $i \in I$ will not intersect
transversally at $X \in \T_N$.
\end{proof}

\section{Combinatorial
patchworking}\label{S:tcon}

The {\it combinatorial patchworking}, also called {\it
$T$-construc\-tion}, is a particular case of the Viro method. The
general Viro method starts with a {\it convex} polyhedral subdivision
of a polytope $\Delta$ contained in the positive orthant $(\RR_+)^n$
of $\RR^n$.  Recall that in this paper all polytopes, including those
of a polyhedral subdivision, are lattice polytopes.  Here, the ambient
lattice is $\Z^n$.

\begin{definition}
A polyhedral subdivision of a polytope $\Delta$ of dimension $n$ is
called {\it convex} (or {\it coherent}) if there exists a convex
piecewise-linear function $\nu : \Delta \to \RR$ whose maximal domains
of linearity coincide with the $n$-polytopes of the subdivision.
\end{definition}

We begin with a brief description of the combinatorial patchworking in
the hypersurface case (see, for example, \cite{IteVir}, \cite{vir6} or
\cite{GKZ}).  Let $\Delta \subset (\RR_+)^n$ be a polytope of maximal
dimension $n$.  Start with a convex triangulation $\calS$ of $\Delta$
and a sign distribution $\delta: \mbox{vert}(\calS) \rightarrow \{\pm
1\}$ at the vertices of $\calS$.  Let $\nu:\Delta \rightarrow \R$ be
any function which certifies the convexity of $\calS$ 
and consider the polynomial
$$f_t(x)= \sum_{\mbox{vert}(\calS)} \delta(w) t^{\nu(w)} x^w$$
where the sum is taken over the set of vertices of $\calS$.
Such a polynomial is called a {\it T-polynomial}.

Denote by $s_{(i)}$ the reflection about the $i$-th coordinate
hyperplane in $\RR^n$.  Let $\Delta^*$ be the union of the $2^n$
symmetric copies of $\Delta$ via compositions of these reflections and
extend $\calS$ uniquely to a triangulation $\calS^*$ which is
symmetric with respect to the coordinate hyperplanes.  Extend the sign
distribution $\delta$ to a sign distribution $\delta^*$ at the
vertices of $\calS^*$ so that a vertex of $\calS^*$ and its image
under a reflection $s_{(i)}$ have the same sign if and only if the
$i$-th coordinate of the vertex is even.  If $\sigma$ is an
$n$-simplex of $\calS^*$ whose vertices have different signs, select
the hyperplane piece which is the convex hull of the middle points of
the edges of $\sigma$ with endpoints of opposite signs.  The union of
all these selected pieces produces a piecewise-linear hypersurface
$H^*$ in $\Delta^*$.

We perform identifications on the boundary of $\Delta^*$ in the
following way.  Let $\Gamma$ be any proper face of $\Delta$ and
consider the cone generated by all outward real vectors which are
orthogonal to the facets of $\Delta$ incident to $\Gamma$.  The
integer vectors in this cone form a finitely generated
semi-goup. Identify two points lying on two symmetric copies of
$\Gamma$ whenever they are symmetric via ${s_{(1)}}^{v_1} \circ
{s_{(2)}}^{v_2} \circ \dots \circ {s_{(n)}}^{v_n}$ for some $v = (v_1,
\dots, v_n)$ in this semi-group. Denote by $ \widetilde{\Delta}$ the
result of these identifications.  By a classical result (see, for
example, \cite{GKZ} Theorem~5.4 p. 383 \cite{Stuint} Proposition 2),
there is an homeomorphism between the real part $\RR X_\Delta$ of
$X_\Delta$ and $\widetilde{\Delta}$.  Moreover, this homeomorphism can
be choosen so that it respects the stratification by torus orbits in
the sense that the real torus orbit corresponding to a
face $\Gamma$ of $\Delta$ is sent to the image under the previous
identifications of the union of the symmetric copies of the interior
of $\Gamma$. In particular, the dense real torus $(\R^*)^n \subset \RR
X_\Delta$ is sent to the union of the symmetric copies of the interior
of $\Delta$.  Denote by $\widetilde{H}$ the image of $H^*$ in
$\widetilde{\Delta}$.

\begin{thm}[T-construction, O. Viro]\label{Tcons}
For $t>0$ sufficiently small, the polynomial $f_t$ is nondegenerate.
Moreover, there exists an homeomorphism $\RR X_{\Delta} \rightarrow \widetilde{\Delta}$ which respects the stratification by torus orbits
and induces an homeomorphism between the real part of the hypersurface $\bar{Z}_f \subset X(\Delta)$ and $\widetilde{H}$.
\end{thm}

We now describe the extension of the combinatorial patchworking to the
case of complete intersections due to B. Sturmfels~\cite{Stu}.  Start
with $k \geq 2$ polytopes $\Delta_1, \dotsc,\Delta_k$ in $(\RR_+)^n$.
Assume that each $\Delta_i$ comes with a convex polyhedral subdivision
$\calS_i$ induced by a convex piecewise-linear map $\nu_i: \Delta_i
\rightarrow \R$.  These functions $\nu_1 ,\dotsc ,\nu_k$ define a
convex polyhedral subdivision of the Minkowsky sum
$\Delta=\Delta_1+\cdots+\Delta_k$ in the following way (see
\cite{Stuint}, \cite{Sturesultant} or \cite{Bi1}). Let
$\hat{\Delta}_i$ be the convex hull of the set $\{(x,\nu_i(x)), x \in
\Delta_i\}$ in $\RR^n \times \RR$ . Let $\hat{\Delta} \subset \RR^n
\times \RR$ be the Minkowski sum ${\hat{\Delta}}_1 + \cdots +
\hat{\Delta}_k$. Let $\calM\calS$ be the convex polyhedral subdivision
of $\Delta$ induced by $\nu$. Each lower face $\hat{\Gamma}$ of
$\hat{\Delta}$ can be uniquely written as a Minkowsky sum
$\hat{\Gamma}_1+\cdots+\hat{\Gamma}_k$ of lower faces of
$\hat{\Delta}_1,\dotsc, \hat{\Delta}_k$. Projecting to $\Delta$, this
gives a representation of each polytope $\Gamma$ of $\calM\calS$ as
$\Gamma = \Gamma_1+ \cdots + \Gamma_k$ with $\Gamma_i \in \calS_i$ for
$i=1,\dotsc,k$.  Such a representation is not unique in general, and
we shall always use the one obtained by projecting lower faces of
$\hat{\Delta}$.  The polyhedral subdision $\calM\calS$ together with
the associated representation of each of its polytopes is called a
convex or coherent {\it mixed subdivision}.  Sturmfels' theorem
requires the following genericity condition.  Namely, assume that each
subdivision $\calS_i$ is a triangulation and that
$$\dim \Gamma = \dim \Gamma_1 + \cdots+ \dim \Gamma_k$$ for any
$\Gamma = \Gamma_1 + \cdots + \Gamma_k \in \calM\calS$ with $\Gamma_i
\in \calS_i$.  We call such a mixed subdivision a convex { \it tight
mixed subdivision}.  (See~\cite{Bi1, Bi2} for other versions of the
Viro method for complete intersections).

Suppose now that for $i=1,\dotsc,k$ a sign distribution $\delta_i:
\mbox{vert}(\calS_i) \rightarrow \pm 1$ is given.  Consider the
T-polynomials associated with these data

$$f_{i,t}(x)=\sum_{\mbox{vert}(\calS_i)} \delta_i(w) t^{\nu_i(w)} x^w \, .$$

Extend $\calM\calS$ to a subdivision $\calM\calS^*$ of $\Delta^*$ by
means of the reflections about coordinate hyperplanes. Hence
$\calM\calS^*$ consists of the polytopes
$s(\Gamma)=s(\Gamma_1)+\cdots+s(\Gamma_k)$ where $s$ is a composition
of coordinate hyperplane reflections and $\Gamma=
\Gamma_1+\cdots+\Gamma_k \in \calM\calS$ ($\Gamma_i \in \calS_i$).
Extend $\delta_i$ to a sign distribution $\delta_i^*$ at the vertices
of $\calS_i^*$ using the rule described above.  Define a sign
distribution $\delta: \mbox{vert}(\calM\calS) \rightarrow \{\pm 1\}^k$
by assigning $(\delta_1(v_1), \dotsc,\delta_k(v_k))$ to each vertex
$v$ of $\calM\calS$ with representation $v=v_1+\cdots+v_k$. Extend
$\delta$ to a sign distribution $\delta^*$ at the vertices of
$\calM\calS^*$ so that the $i$-th sign of a symmetric copy $s(v)$ of
$v=v_1+\dotsc+v_k$ is $\delta_i^*(s(v_i))$.

For any $i=1,2,\dotsc,k$, let $H_i^* \subset \Delta_i^*$ be the piecewise-linear hypersurface
constructed via the combinatorial patchworking from $\calS_i$ and $\delta_i$. Let $H_i^{\Delta, *} \subset \Delta^*$
be the union over all polytopes $s(\Gamma)=s(\Gamma_1)+\dotsc+s(\Gamma_k) \in \calM\calS^*$
of $\oplus_{j \neq i} s(\Gamma_j) \, +H_i^* \cap s(\Gamma_i) \subset s(\Gamma)$.
Let $\widetilde{H_i^{\Delta}}$ denote the image of $H_i^{\Delta}$ in $\widetilde{\Delta}$.

\begin{thm}[B. Sturmfels]\label{St}
For $t>0$ sufficiently small the collection $(f_{1,t},\dotsc,f_{k,t})$ is nondegenerate.
Moreover, there exists an homeomorphism $\RR X_{\Delta} \rightarrow \widetilde{\Delta}$ which respects the stratification by torus orbits
and induces for each $i$ an homeomorphism between the real part of the hypersurface $\bar{Z}_f^{\Delta} \subset X(\Delta)$ and $\widetilde{H_i^{\Delta}}$.
\end{thm}

These two versions -- for hypersurfaces and complete intersections -- of the combinatorial patchworking
are related by
the so-called {\it combinatorial Cayley trick}.
Consider the Cayley polynomial $F_t \in \R[x,y]$ associated with $(f_{1,t},\dotsc,f_{k,t})$
$$
F_t(x,y)=\sum_{i=1}^k y_i \, f_{i,t}(x) \,.
$$ Its Newton polytope is the Cayley polytope $C(\Delta_1, \ldots ,
\Delta_k) \subset \R_+^{n+k}$.  Let $(a,b)$ be coordinates on
$\R^{n+k}=\R^n \times \R^k$. Consider the subspace $B$ of $\R^{n+k}$
defined by $b_1=b_2=\cdots=b_k=1/k$ and identify it with $\R^n$ via
the projection $(a,b) \mapsto a$.  This identifies $B \cap C(\Delta_1,
\ldots , \Delta_k)$ with $\Delta=\Delta_1+\cdots+\Delta_k$ dilated by
$1/k$.  Note that the space defined by $b_i=1$ and $b_j=0$ for $j \neq
i$ intersects $C(\Delta_1, \ldots , \Delta_k)$ along a face which can
be identified with $\Delta_i$ via the projection.  Consider a
polyhedral subdivision of $C(\Delta_1,\ldots , \Delta_k)$.  If $F$ is
a polytope of maximal dimension $\dim \Delta +k-1$ in this
subdivision, then it intersects the space defined by $b_i=1$ and
$b_j=0$ for $j \neq i$ along a nonempty face $F_i$, which projects to
a (nonempty) subpolytope $\Gamma_i$ of $\Delta_i$.  Then $F \cap B $
is identified via the projection with the polytope
$\Gamma=\Gamma_1+\cdots+\Gamma_k \subset \Delta$ dilated by $1/k$.
This gives a correspondence between polyhedral subdivisions of
$C(\Delta_1, \ldots , \Delta_k)$ and mixed subdivisions of
$\Delta=\Delta_1+\cdots+\Delta_k$. It is easily seen that
triangulations of $C(\Delta_1, \ldots , \Delta_k)$ are sent to tight
mixed subdivision via this correspondence.  The following result can
be found, for example, in \cite{Sturesultant}.

\begin{prop}\label{P:Cayley-trick}
The correspondence described above is a bijection between the set of
convex polyhedral subdivision of $C(\Delta_1, \ldots , \Delta_k)$ and
the set of mixed subdivisions of
$\Delta=\Delta_1+\cdots+\Delta_k$. Precisely, let
$\nu: C(\Delta_1, \ldots , \Delta_k) \rightarrow \R$ be
any convex piecewise-linear function and let $\nu_i$
denote its restriction to $\Delta_i$ identified with a face of
$C(\Delta_1, \ldots , \Delta_k)$ via the projection $(a,b) \mapsto a$.
Then the correspondence described above sends the coherent polyhedral subdivision of $C(\Delta_1,
\ldots , \Delta_k)$ defined by $\nu$
to the coherent mixed subdivision of $\Delta$ defined by $(\nu_1,\dotsc,\nu_k)$.
\end{prop}

Note that in the situation of Theorem~\ref{St}, the Cayley polynomial
$F_t$ is a T-polynomial.  The non degeneracy of
$(f_{1,t},\dotsc,f_{k,t})$ in Theorem~\ref{St} follows from
Proposition~\ref{P:Cayleytrick} and Theorem~\ref{Tcons} applied to
$F_t$.

\section{Standard definitions and properties in tropical geometry}
\label{S:standard}
The setting and notation here are the same as in \cite{Be1}.  A
detailed exposition can be found in~\cite{MCounting} and in
~\cite{IMS}, for example.  Let $\KK$ be the field of Puiseux series.
An element of $\KK$ is a series $g(t)=\sum_{r\in R} b_r t^r$ where
each $b_r$ is a complex number and $R \subset \mathbb{Q}$ is bounded
from below and contained in an arithmetic sequence.  Consider the
valuation $\val (g(t)):=\min \{r \, | \, b_r \neq 0 \}$. Using
Mikhalkin's conventions, we rather use minus the valuation
$v(g):=-\val(g)$. Define
\[
\begin{array}{rrrl} 
 V: &(\KK ^*)^n &\longrightarrow &\RR^n \\
    & z  &  \longmapsto &(v(z_1), \dots, v(z_n)).
\end{array}
\]

Let $f$ be a polynomial in $\KK[z_1,\dotsc,z_n] =\KK[z]$.  It is of the form
$f(z)=\sum_{\omega \in A} \;c_\omega z^\omega$ with $A$ a finite
subset of $\ZZ^n$ and $c_\omega \in \KK^*$. Let $Z_f=\{z\in (\KK
^*)^n \, | \, f(z)=0 \}$ be the zero set of $f$ in $(\KK ^*)^n$.

\begin{definition}
  The tropical hypersurface $Z_f^{\trop}$ associated to $f$ is the
  closure (in the usual topology) of the image under $V$ of $Z_f$:
$$ Z^{\trop}_f = \overline{V(Z_f)} \subset \RR^n.$$
\end{definition}
There are other equivalent definitions of a tropical hypersurface.
Namely, define
\[
\begin{array}{rrrl} 
 \nu: & A &\longrightarrow &\RR \\
                & \omega  &  \longmapsto & - v (c_\omega)
\end{array}
\] 
Its Legendre transform is the
piecewise-linear convex function
\[
\begin{array}{rrrl} 
\mathcal{L}(\nu): &\RR^n &\longrightarrow &\RR \\
                & x  &  \longmapsto & \max_{\omega\in A} (x \cdot \omega - \nu(\omega)) 
\end{array}
\]

\begin{thm}[Kapranov]\label{T:Kapranov}
The tropical hypersurface $Z^{\trop}_f$ is the corner locus of  $\mathcal{L}(\nu)$.
\end{thm} 
\noindent The corner locus of $\mathcal{L}(\nu)$ is the set of points at which it
is not differentiable. Another way to define a tropical hypersurface is to use the
{\it tropical semiring} $\Rtrop$, which is $\R \cup \{-\infty\}$ endowed
with the following tropical operations. The tropical addition
of two numbers is the maximum of them, and thus its neutral element is $-\infty$.
The tropical multiplication is the ordinary addition with the convention that $x+(-\infty)=-\infty+x=-\infty$.
Removing the neutral element for the tropical addition,
we get the {\it one dimensional tropical torus} $\Trop:=\R=\Rtrop \setminus \{-\infty\}$.
A multivariate tropical polynomial is a polynomial in
$\R[x_1,\dotsc,x_n]$ where the addition and multiplication are the tropical ones (strictly speaking, the
coefficients are in $\Rtrop$, but as usual we omit the monomials whose
coefficients are the neutral element for the addition). Hence, a tropical polynomial is given by a
maximum of finitely many affine functions whose linear parts have
integer coefficients and constant parts are real numbers.
The tropicalization of a polynomial
$$f(z)=\sum_{\omega \in A} \;c_\omega z^\omega \in \KK[z]$$
where the coefficients $c_\omega \in \KK$ are all nonzero
is the tropical polynomial
$$\mbox{Trop}(f)(z)=\sum_{\omega \in A} \; v(c_\omega) z^\omega \in \R[z].$$
This tropical polynomial coincides with the piecewise-linear convex function
$\mathcal{L}(\nu)$ defined above.  Therefore, Theorem~\ref{T:Kapranov}
asserts that $Z^{\trop}_f$ is the corner locus of $\mbox{Trop}(f)$.
Conversely, the corner locus of any tropical polynomial is a tropical
hypersurface (just take a polynomial in $\KK[z]$ whose coefficients
have the right valuations). For these reasons, we will sometimes speak
about the tropical hypersurface defined by a polynomial $f$ without
specifying if $f$ is in $\KK[z]$ or if $f$ is a tropical polynomial
(the tropicalization of the latter).

The Newton polytope of the tropical hypersurface $Z^{\trop}_f$ is the
convex hull of $A$ and will be denoted by $\Delta$.  One can associate
to $Z^{\trop}_f$ a polyhedral subdivision $\calS$ of $\Delta$ in the
following way. Let $ \hat{\Delta}\subset \RR^n\times\RR$ be the convex hull of all points
$(\omega,v(c_\omega))$ with $\omega \in A$. Define

\begin{equation}\label{E:dualsubdivision}
\begin{array}{rrrl} 
\hat{\nu} : & \Delta &\longrightarrow &\RR \\
                & x &  \longmapsto & \min \{y \, | \,  (x,y) \in \hat{\Delta}\}.
\end{array}
\end{equation}

The domains of linearity of $\hat{\nu}$ form a convex polyhedral
subdivision $\mathcal S$ of $\Delta$. The hypersurface $Z_f^{\trop}$
is an $(n-1)$-dimensional piecewise-linear complex which induces a
polyhedral subdivision $\Xi$ of $\RR^n$. We will call {\it cells} the
elements of $\Xi$.  Note that these cells have rational slopes.  The
$n$-dimensional cells of $\Xi$ are the closures of the connected
components of the complement of $Z_f^{\trop}$. The lower dimensional cells of $\Xi$ are
contained in $Z_f^{\trop}$ and we will just say that they are cells of
$Z_f^{\trop}$.  Both subdivisions $\calS$ and $\Xi$ are dual in the
following sense.  There is a one-to-one correspondence between $\Xi$
and $\calS$, which reverses the inclusion relations, and such that if $\sigma \in \calS$
corresponds to $\xi \in \Xi$ then
\begin{enumerate}
\item $\dim \xi +\dim \sigma=n$,
\item the cell $\xi$ and the polytope $\sigma$ span orthonogonal real
affine spaces,
\item the cell $\xi$ is unbounded if and only if $\sigma$ lies on a
proper face of $\Delta$.
\end{enumerate}
Note that under this correspondence the cells of $Z_f^{\trop}$
correspond to positive dimensional polytopes of $\calS$.  We now
underline some similarities between complex toric hypersurfaces and
tropical hypersurfaces.  As in the complex case, we can start with a
polynomial $f$ whose exponent vectors belong to a lattice $M \simeq
{\Z}^n$.  Then, in view of the definition of $\mathcal{L}(\nu)$ and
Theorem~\ref{T:Kapranov}, the tropical hypersurface lies in the real
vector space $N_{\R} \simeq \R^n$ generated by the lattice $N$ dual to
$M$. This real vector space $N_{\R}$ can be interpreted as the
tropical torus $\Trop_N$ associated with the lattice $N$, so that
$Z_f^{\trop} \subset \Trop_N=N_{\R}$ is in fact a toric tropical
hypersurface. The polynomial $f$ also defines a toric tropical
hypersurface in $N(\Delta)_{\R} \simeq \R^{\dim \Delta}$ and
$Z_f^{\trop} \subset N_{\R}$ is the product of this hypersurface with
the tropical torus $\simeq \R^{n-\dim \Delta}$ associated with (the
dual of) a complement of $M(\Delta)$ in $M$.  The unbounded cells of
$Z_f^{\trop}$ gives rise to toric tropical hypersurfaces defined by
truncations of $f$ to faces of $\Delta$.  Namely, consider a face
$\Gamma$ of $\Delta$ and let $\gamma \subset N$ be the semigroup
formed by all elements of $N$ which are identically zero on
$M(\Gamma)$ and are negative on any vector $w=m'-m \in M(\Delta)$ with
$m' \in \Delta \setminus \Gamma$ and $m \in \Gamma$ (in other words,
$\gamma$ consists of all integer vectors of $N$ orthogonal to $\Gamma$
and going outside $\Delta$).  Note that $N(\Gamma)$ is the quotient
$\frac{N}{\gamma + (-\gamma)}$, where $\gamma + (-\gamma)$ is the
subgroup of $N$ generated by the semigroup $\gamma$.  Consider the
unbounded cells of $Z_f^{\trop}$ which intersect any hyperplane $\{w
\in M_{\R} \, | \, v \cdot w =c\}$ with $c$ big enough and $v$ in
$\gamma$. The cells of the tropical hypersurface
$Z_{f^{\Gamma}}^{\trop} \subset N(\Gamma)_{\R}$ are exactly the images
of these cells under the quotient map $N_{\R} \rightarrow
N(\Gamma)_{\R}$. Comparing with the classical complex situation, this
leads to the notion of {\it tropical variety} $\Trop_{\Delta}$
associated with $\Delta$ with properties analogous to those of the
complex projective toric variety $X_{\Delta}$. Geometrically, one can
think about $\Trop_{\Delta}$ as being the image of $\Delta$ by the
composition of a translation and a dilatation, so that
$\bar{Z}_{f}^{\trop} \subset \Trop_{\Delta}$ can be obtained from
$Z_{f}^{\trop}$ by cutting the the unbounded cells of
${Z}_{f}^{\trop}$ along the faces of $\Trop_{\Delta}$.
For the sake of completness, we recall the definition of a tropical variety in $N_{\R}$
(see, for example, \cite{EKL04} and \cite{Gathmann}).

\begin{definition}
A tropical variety in $N_{\R}$ is the closure of the image under $V$
of the zero set of an ideal $I \subset \KK[z_1,\dotsc,z_n] =\KK[z]$.
We will denote this tropical variety by $Z^{\trop}_I$.
\end{definition}

It turns out that the tropical variety $Z^{\trop}_I$ is the common
intersection of all tropical hypersurfaces $Z^{\trop}_f$ for $f \in I$
(see~\cite{SS,P}). There exists a finite number of polynomials $f_1,\dotsc,f_k \in
I\subset \KK[z]$ so that $Z^{\trop}_I$ is the common intersection of
the corresponding tropical hypersurfaces (see~\cite{TADBR, HT}).  Such a
collection of polynomials is called a {\it tropical basis} of
$Z^{\trop}_I$.  On the other hand, it is known that the common
intersection of tropical hypersurfaces is not always a tropical
variety.

\section{Intersection multiplicity numbers between tropical hypersurfaces}
\label{S:mult}

Recall that all polytopes under consideration have vertices in the
underlying lattice $M \simeq \Z^n$.  A $k$-dimensional simplex
$\sigma$ with vertices $m_0, m_1 ,\dotsc,m_k$ is called {\it
primitive} if the vectors $m_1-m_0,\dotsc,m_k-m_0$ form a basis of the
lattice $M(\sigma)$, or equivalently, if these vectors can be
completed to form a basis of $M$. Obviously, the faces of a primitive
simplex are themselves primitives simplices.

Consider a $k$-dimensional vector subspace of $M_{\R}$ with rational
slopes.  It intersects $M$ in a saturated subgroup $\gamma$ of rank
$k$ and coincides with the real vector space $\gamma_{\R}$ generated
by $\gamma$.  Any basis of $\gamma$ produces an isomorphism between
$\gamma$ and $\Z^k$, and then by extension an isomorphism between
$\gamma_{\R}$ and $\R^k$.  Let $\mbox{Vol}_{\gamma}$ be the volume
form on $\gamma_{\R}$ obtained as the pull-back via such an
isomorphism of the usual Euclidian $k$-volume on $\R^k$. For
simplicity, we will write $\mbox{Vol}_k$ instead of
$\mbox{Vol}_{\gamma}$ since the lattice $\gamma$ will be clear from
the context.  Note that $\mbox{Vol}_k$ does not depend on the
isomorphism $\gamma \simeq \Z^k$ since two basis of $\gamma$ are
obtained from each other by integer invertible linear map which has
determinant $\pm 1$.  Any basis $(\gamma_1,\dotsc,\gamma_k)$ of
$\gamma$ generate a $k$-dimensional parallelotope $P \subset
\gamma_{\R}$ (isomorphic to the cube $[0,1]^k \subset \R^k$) called
{\it fundamental parallelotope} of $\gamma$ and which verifies
$\mbox{Vol}_k (P)=1$. Two primitive $k$-simplices on $\gamma_{\R}$
have the same volume under $\mbox{Vol}_k$ (they are interchanged by an
invertible integer linear map), and this volume is $\frac{1}{k!}$
since a fundamental parallelotope of $\gamma$ can be subdivided into
$k!$ primitive $k$-simplices.  We will often use the normalized volume
$$\mbox{vol}_k (\, \cdot \, ) := k ! \cdot \mbox{Vol}_k(\, \cdot \,
)$$ on $\gamma_{\R}$. This normalized volume takes all nonnegative
integer values on polytopes (with vertices in $\gamma$), and we have
$\mbox{vol}_k(\sigma)=1$ for a polytope $\sigma$ if and only if
$\sigma$ is a $k$-dimensional primitive simplex.  We will use the
following elementary fact.

\begin{rem}\label{R:indexdet}
Let $\gamma$ be a subgroup of a free abelian group $\Lambda$ of finite rank.
Assume that $\Lambda$ and $\gamma$ have the same rank $k$, so that the index
$[\Lambda: \gamma]$ of $\gamma$ in $\Lambda$ is well-defined.  Then,
for any basis $(\gamma_1,\dotsc,\gamma_k)$ of $\gamma$ and any basis
$e=(e_1,\dotsc,e_k)$ of $\Lambda$ we have
$$[\Lambda: \gamma]=\mbox{Vol}_k(G)=\mbox{vol}_k(g)=|\det(G_{ij})|,$$
where $G$ (resp., $g$) is the $k$-dimensional parallelotope (resp.,
$k$-dimensional simplex) generated by $\gamma_1,\dotsc,\gamma_k$ and
$(G_{ij})$ is the $k \times k$-matrix whose $j$-th column is the
vector of coordinates of $\gamma_j$ with respect to
$(e_1,\dotsc,e_k)$.
\end{rem}

Consider now tropical polynomials $f_1,\dotsc,f_k$ in
$\R[x_1,\dotsc,x_n]$ or more generally in $\R[M]$ with $M \simeq
\Z^n$. Denote by $\Delta_i$ the Newton polytope of $f_i$.  Recall that
each tropical hypersurface $Z_{f_i}^{\trop}$ defines a piecewise
linear polyhedral subdivision $\Xi_i$ of $N_{\R}$ which is dual to a
convex polyhedral subdivision $\calS_i$ of $\Delta_i$.  The union of
these tropical hypersurfaces defines a piecewise-linear polyhedral
subdivision $\Xi$ of $N_{\R}$.
Any non-empty cell of $\Xi$ can be written as
\begin{equation}\label{E:cell}
\xi=\bigcap_{i=1}^k \xi_i
\end{equation}
with $\xi_i \in \Xi_i$ for $i=1,\dotsc,k$. Any cell $\xi \in \Xi$ can
be uniquely written in this way if one requires that $\xi$ lies {\it
in the relative interior} of each $\xi_i$. We shall always refer to this unique
form. Denote by $\calM \calS$ the mixed subdivision of
$\Delta=\Delta_1+\cdots+\Delta_k$ induced by the tropical polynomials
$f_1,\dotsc,f_k$. Recall that any polytope $\sigma \in \calM \calS$
comes with a privileged representation
$$\sigma=\sigma_1+\cdots+\sigma_k$$ 
with $\sigma_i \in \calS_i$.  The
above duality-correspondence applied to the (tropical) product of the
tropical polynomials gives rise to the following fact.

\begin{prop}\label{P:nonempty}
There is a one-to-one duality correspondence between $\Xi$ and
$\calS$, which reverses the inclusion relations, and such that if
$\sigma \in \calM\calS$ corresponds to $\xi \in \Xi$ then
\begin{enumerate}
\item if $\xi=\bigcap_{i=1}^k \xi_i$
with $\xi_i \in \Xi_i$ (and $\xi$ lies in the relative interior of $\xi_i$) for $i=1,\dotsc,k$, then $\sigma$ has
representation $\sigma=\sigma_1+\cdots+\sigma_k$ where each $\sigma_i$
is the polytope dual to $\xi_i$.

\item $\dim \xi +\dim \sigma=n$,

\item the cell $\xi$ and the polytope $\sigma$ span orthonogonal real
affine spaces,

\item the cell $\xi$ is unbounded if and only if $\sigma$ lies on a
proper face of $\Delta$.
\end{enumerate}

\end{prop}

We put weights on the cells of each subdivision $\Xi_i$ in the
following way.  If $\xi_i \in \Xi_i$ is a cell of maximal dimension
$n$ (which means that $\xi_i$ is not a cell of the tropical
hypersurface $Z_{f_i}^{\trop}$), then its weight is defined by
$w(\xi_i):=0$.  If $\xi_i \in \Xi_i$ is a cell of positive codimension
$d_i$, then
$$
w(\xi_i):=\mbox{vol}_{d_i} (\sigma_i)$$
where $\sigma_i \in \calS_i$ is the polytope corresponding to $\xi_i$.
We now define weights on the cells of $ \Xi$ in the following way.
Consider a cell $\xi \in \Xi$
$$
\xi=\bigcap_{i=1}^k \xi_i
$$ where $\xi_i \in \Xi_i$ for $i=1,\dotsc,k$ (and $\xi$ lies in the
relative interior of each $\xi_i$). Let $\sigma_i \in \calS_i$ be the
polytope corresponding to $\xi_i$. Set $d_i:=\codim \xi_i=\dim
\sigma_i$ and $d:=\codim \xi=\dim \sigma$.  Recall that for a polytope
$P \subset M_{\R}$, we denote by $M(P)$ the subgroup of $M$ consisting
of all integer vectors which are parallel to $P$.

\begin{definition}\label{D:weights}
The weight of $\xi$ is defined as follows.

\begin{itemize}
\item ({\it Tranversal case.}) If $d_1+\cdots+d_k=d$, then
$$
\begin{array}{lll}
w(\xi) & = & \left(\prod_{i=1}^k w(\xi_{i}) \right) \cdot [M(\sigma) :
 M(\sigma_1)+\cdots+M(\sigma_k)] \\
 & & \\
 & = & \left(\prod_{i=1}^k
 \mbox{vol}_{d_i}(\sigma_{i}) \right) \cdot [M(\sigma)
 :M(\sigma_1)+\cdots+M(\sigma_k)] 
\end{array}
$$

\item ({\it General case.})  Translate the tropical hypersurfaces by
small generic vectors so that all intersections emerging from $\xi$
are transversal intersections. Define $w(\xi)$ as the sum of the
weights at the transversal intersections emerging from $\xi$ and which
are cells of codimension $d$.
\end{itemize}
\end{definition}

Our weights are similar to those introduced by Mikhalkin
in~\cite{Mtropappli} in order to define tropical cycles. Note
that in~\cite{Mtropappli} only top-dimensional cells are equipped with
weights. In our situation, the top-dimensional cells of the cycle
corresponding to the intersection of our tropical hypersurfaces are
cells $\xi \in \Xi$ of codimension $d=k$. It follows straightforwardly
from the definitions and Lemma~\ref{L:index} below that on these
top-dimensional cells our weights coincide with those of Mikhalkin.
We will show in Theorem~\ref{T:weights} that our weight does not
depend (in the non transversal case) on the translation vectors. It is
then natural to interpret $w(\xi)$ as being the {\it intersection
multiplicity number} between the tropical hypersurfaces
$Z_{f_1}^{\trop}, \dotsc, Z_{f_k}^{\trop}$ along the cell $\xi$.

\begin{lemma}\label{L:index}
Let $\gamma_1$ and $\gamma_2$ be saturated subgroups of a free group
$N$ such that $\gamma_1+\gamma_2$ and $N$ have same rank.  Then the
index of $\gamma_1+\gamma_2$ in $N$ satisfies to
$$[N : \gamma_1+\gamma_2] =[ (\gamma_1 \cap \gamma_2)^{\bot}:
\gamma_1^{\bot}+\gamma_2^{\bot}],$$ where $\gamma^{\bot}$ denotes the
subgroup of the dual lattice $M=\Hom_{\Z}(N,\Z)$ consisting of all
elements of $M$ which vanish on a subgroup $\gamma$ of $N$.
\end{lemma}
\begin{proof}
If $\gamma_1 \cap \gamma_2=\{0\}$ then $(\gamma_1 \cap
\gamma_2)^{\bot}=M$ and the corresponding equality has been proven
in~\cite{Katz2}.  The general case reduces to this case in the
following way.  Let $n$ be the rank of $N$. If $\gamma_1$ and
$\gamma_2$ are saturated then so is $\gamma_1 \cap \gamma_2$. This
implies that the quotient $N/(\gamma_1 \cap \gamma_2)$ is a free group
of rank $n-\rk(\gamma_1 \cap \gamma_2)$.  We have a group isomorphism
$$
\frac{N}{\gamma_1+\gamma_2} \simeq
\frac
{N /(\gamma_1 \cap \gamma_2)}
{ (\gamma_1+\gamma_2) / (\gamma_1 \cap \gamma_2)}$$
and also
$$(\gamma_1+\gamma_2) / (\gamma_1 \cap
\gamma_2)=\frac{\gamma_1}{\gamma_1 \cap
\gamma_2}+\frac{\gamma_2}{\gamma_1 \cap \gamma_2}.$$ The group dual to
$N/(\gamma_1 \cap \gamma_2)$ is isomorphic to $(\gamma_1 \cap
\gamma_2)^{\bot} \subset M$.  It remains to note that if $\gamma_1$
and $\gamma_2$ are saturated subgroups of $N$, then for $i=1, 2$ the
subgroup $\frac{\gamma_i}{\gamma_1 \cap \gamma_2}$ of
$\frac{N}{\gamma_1+\gamma_2}$ is also saturated.
\end{proof}

Let $P_1,\dotsc,P_{\ell}$ be polytopes with vertices in a saturated
lattice $\gamma$ of rank $\ell$.  The map
$(\lambda_1,\dotsc,\lambda_{\ell}) \mapsto
\mbox{Vol}_{\ell}(\lambda_1P_1+\cdots+\lambda_{\ell}P_{\ell})$ is a
homogeneous polynomial map of degree $\ell$.
The coefficient of the monomial $\lambda_1 \cdots \lambda_{\ell}$ 
is called the {\it mixed volume} of $P_1,\dotsc,P_{\ell}$ and is denoted by
$$MV_{\ell}(P_1,\dotsc,P_\ell).$$ A famous theorem due to Bernstein
states that this mixed volume is the number of solutions in the torus
associated with the lattice $\gamma$ of a generic polynomial system
$f_1=\dotsc=f_{\ell}=0$ where each $f_i$ has $P_i$ as Newton
polytope. Note that $MV_{\ell}(P_1,\dotsc,P_\ell)=0$ if
$P=P_1+\cdots+P_{\ell}$ has not full dimension $\ell$ or if at least
one $P_i$ has dimension zero. We may also consider mixed volumes
associated with any number $m \leq \ell$ of polytopes among
$P_1,\dotsc,P_{\ell}$ (see~\cite{USA}). Namely, let $P_1,\dotsc,P_{m}$
be $m \leq \ell$ polytopes with vertices in a lattice of rank $\ell$
and let $\underline{t}=(t_1,\dotsc,t_m)$ be a collection of positive integer numbers such that
$\sum_{i=1}^m t_i={\ell}$. Then define

$$MV_{\ell}(P_1,\dotsc,P_m \, ; \, \underline{t}) :=
MV_{\ell}(\underbrace{P_1,\dotsc,P_1}_{t_1}, \dotsc \dotsc,
\underbrace{P_m,\dotsc,P_m}_{t_m}),$$ where on the right each $P_i$ is
repeated $t_i$ times. We note that $\frac{1}{t_1!\cdots t_m!} \cdot MV_{\ell}(P_1,\dotsc,P_m \, ; \, \underline{t})$
is the coefficient of $\lambda_1^{t_1}\cdots\lambda_m^{t_m}$ in the homogeneous
degree $\ell$ polynomial map $(\lambda_1,\dotsc,\lambda_m) \mapsto \mbox{Vol}_{\ell}(\lambda_1P_1+\cdots+\lambda_mP_m)$ (see \cite{USA},
page 327).

We are now able to state a formula for the weight $w(\xi)$ defined
above which shows in particular that this weight does not depend on
the chosen translation vector.

\begin{thm}\label{T:weights}
Let $\xi$ be a cell of $\Xi$ and $w(\xi)$ be its weight as defined in
Definition~\ref{D:weights}.  If $\xi$ is not a cell of $\cap_{i=1}^k
Z_{f_i}^{\trop}$, or equivalently if at least one $d_i$ is zero, then
$w(\xi)=0$. Assume now that all $d_i$ are positive integer numbers.

\begin{itemize}
\item If the tropical hypersurfaces intersect transversally along
$\xi$, which means that $d=d_1+\cdots+d_k$, then letting
$\underline{d}:=(d_1,\dotsc,d_k)$ we have
\smallskip

\begin{equation}\label{E:transversal}
w(\xi) \;  =  \; MV_d(\sigma_1,\dotsc,\sigma_k ; \underline{d})
\end{equation}
\smallskip

\item In the general case, we have $d \leq d_1+\cdots+d_k$ and

\begin{equation}\label{E:nontransversal}
w(\xi)= \sum_{ \underline{t} \, , \, t_1+\cdots+t_k=d}
MV_d(\sigma_1,\dotsc,\sigma_k ; \underline{t})
\end{equation}
\end{itemize}

As a particular case, if $d=\codim \xi=k$ then

\begin{equation}\label{E:nontransversaltop}
w(\xi)=MV_k(\sigma_1,\dotsc,\sigma_k).
\end{equation}
\end{thm}

Note that in the hypersurface case we have $\xi=\xi_1$ and thus
$w(\xi)=MV_{d_1}(\sigma_1,\dotsc,\sigma_1) = d_1!
\cdot \mbox{Vol}_{d_1}(\sigma_1)=\mbox{vol}_{d_1}(\sigma_1)$ as
required.
Before giving a proof of Theorem~\ref{T:weights}, we need
some intermediate results.

A Minkowsky sum $Q_1+\cdots+Q_{\ell}$ of polytopes such that $\dim
(Q_1+\cdots+Q_{\ell})=\dim Q_1+\cdots+\dim Q_{\ell}$ is called {\it a
direct Minkowsky sum} and is denoted by $Q_1 \oplus \cdots \oplus
Q_{\ell}$.  A convex mixed subdivision $\calM \calS$ of a polytope
$P=P_1+\cdots+P_{\ell}$ is called {\it pure} if for any polytope $Q
\in \calM \calS$ with (privileged) representation
$Q=Q_1+\cdots+Q_{\ell}$ we have $Q=Q_1 \oplus \cdots \oplus Q_{\ell}$.
A convex tight mixed subdivision is a convex pure mixed
subdivision with the additional property that each $Q_i$ is a
simplex. Tight and pure convex mixed subdivisions are generic within
all convex mixed subdivisions of a given collection of polytopes.

\begin{lemma}\label{L:mixedvolumesum}
Let $P=P_1+\cdots+P_{\ell} \subset M_{\R} \simeq \R^{\ell}$.
\begin{enumerate}
\item For any convex pure mixed subdivision of $P=P_1+\cdots+P_{\ell}
$, we have
\begin{equation}\label{E:pure}
\mbox{MV}_{\ell}(P_1,\dotsc,P_{\ell})=\sum \mbox{Vol}_{\ell}(Q_1
\oplus \cdots \oplus Q_{\ell})
\end{equation}
where the sum is taken over all polytopes $Q=Q_1 \oplus \cdots \oplus Q_{\ell}$ of the mixed subdivision with
$\dim Q_1= \cdots = \dim Q_{\ell}=1$.
\item More generally, for any convex mixed subdivision of
$P=P_1+\cdots+P_{\ell} $, we have
\begin{equation}\label{E:nonpure}
\mbox{MV}_{\ell}(P_1,\dotsc,P_{\ell})=\sum
\mbox{MV}_{\ell}(Q_1,\dotsc,Q_{\ell})
\end{equation}
where the sum is taken over all polytopes $Q=Q_1+\cdots+Q_{\ell}$ of
the mixed subdivision.
\end{enumerate}
\end{lemma}
\begin{proof}
Formula~\eqref{E:pure} is well-known (see, for example, \cite{USA}, Ch. 7)
and not difficult to prove from the definition of mixed volume given above.
Formula~\eqref{E:nonpure} is a simple consequence
of~\eqref{E:pure}. Indeed, we may perturb slightly functions
$\nu_1,\dotsc,\nu_{\ell}$ determining a given convex mixed subdivision
of $P=P_1+\cdots+P_{\ell} $ so that the new functions induce pure
mixed subdivisions of each polytope $Q=Q_1+\cdots +Q_{\ell}$ of the
initial mixed subdivision. Then these functions define a pure mixed
subdivision of $P=P_1+\cdots+P_{\ell} $ and it remains to
apply~\eqref{E:pure} simultaneously to all these pure mixed
subdivisions.
\end{proof}

\noindent {\it Proof of Theorem~\ref{T:weights}.}  Assume first that
$d=d_1+\cdots+d_k$ (transversal case).  If some $d_i$ is zero, then it
follows directly from Definition~\ref{D:weights} that $w(\xi)=0$.
Assume that $d_i \geq 1$ for $i=1,\dotsc,k$. We prove
Formula~\eqref{E:transversal} with the help of Bernstein's theorem.
Consider a generic polynomial system
\begin{equation}\label{E:system}
f_{1,1}=\cdots=f_{1,d_1}= \cdots \cdots = f_{k,1}=\cdots=f_{k,d_k}=0
\end{equation}
of $d=d_1+\cdots+d_k$ equations where each $f_{i,j}$ has $\sigma_i$ as
Newton polytope.
By Bernstein's theorem, the system~\eqref{E:system} has $MV_d(\sigma_1,\dotsc,\sigma_k ; \underline{d})$
solutions in the complex torus associated with $M(\sigma)$, where
$\underline{d}=(d_1,\dotsc,d_k)$.  Since $\sigma=\sigma_1 \oplus
\cdots \oplus \sigma_k$, the number of solutions to~\eqref{E:system}
in the complex torus associated with $M(\sigma_1)+\cdots+M(\sigma_k)$ is the
product $N=\prod_{i=1}^k N_i$ where $N_i$ is the number of solutions
in the complex torus associated with $M(\sigma_i)$
to the system
$$f_{i,1}=\dotsc=f_{i,d_i}=0.$$  By Bernstein's theorem we have
$N_i=\mbox{MV}_{d_i}(\sigma_i,\dotsc,\sigma_i)= d_i ! \cdot
\mbox{Vol}_{d_i}(\sigma_i)=\mbox{vol}_{d_i}(\sigma_i)$.  Let
$(e_1,\dotsc,e_d)$ be a basis of $M(\sigma)$ and identify the
associated complex torus with $(\C^*)^d$ via this basis.  Let $z_1,\dotsc,z_N$
be the solutions to~\eqref{E:system} in the subtorus of $(\C^*)^d$
associated with $M(\sigma_1)+\cdots+M(\sigma_k)$.  Then the solutions
to~\eqref{E:system} in $(\C^*)^d$ are obtained by solving for each
$z_l$ a system
$$x^{m_i}=z_{l,i} \; , \quad i=1,\dotsc, d,$$ where $m_1,\dotsc,m_d$
are the vectors of coordinates of a basis of
$M(\sigma_1)+\cdots+M(\sigma_k)$ with respect to $(e_1,\dotsc,e_d)$
and $z_l=(z_{l,1},\dotsc,z_{l,d}) \in (\C^*)^d$.  The number of
solutions to such a system is the absolute value of the $(d \times
d)$-determinant $|m_{i,j}|$ which is equal to $[M(\sigma):
M(\sigma_1)+\cdots+M(\sigma)]$. This proves
Formula~\eqref{E:transversal}.

Consider now the general case. We have obviously $d \leq
d_1+\cdots+d_k$.  let $\nu_i :\Delta_i
\rightarrow \R$, $i=1,\dotsc,k$, be the functions given by the tropical hypersurfaces
and which induce the corresponding mixed subdivision $\calM \calS$ of
$\Delta=\Delta_1+\cdots+\Delta_k$.  Denote by $\calS_i$ the convex
polyhedral subdivision of $\Delta_i$ induced by $\nu_i$.  Translations
of the tropical hypersurfaces by a small generic vector correspond to
small perturbations $\tilde{\nu_i}$
of  the functions $\nu_i$ so that for each $i=1,\dotsc,k$ the polyhedral
subdivision of $\Delta_i$ induced by the resulting function
$\tilde{\nu_i}$ coincide with $\calS_i$.  The intersections between
the tropical hypersurfaces which emerge from $\xi$ after such small
perturbations are transversal intersections if and only if the
mixed subdivision of $\sigma=\sigma_1+\cdots+\sigma_k$ induced by
$\tilde{\nu_1}, \dotsc, \tilde{\nu_k}$ is a pure mixed subdivision
$\calM \calS (\sigma)$. Then each polytope $\Gamma \in \calM
\calS(\sigma) $ has a privileged representation

\begin{equation}\label{E:repres}
\Gamma=\Gamma_1 \oplus \cdots \oplus \Gamma_k
\end{equation}
where each $\Gamma_i \in \calS_i$ and the weight of $\xi$ is by
definition the sum of weights of the cells corresponding to polytopes
$\Gamma \in \calM \calS(\sigma)$ such that $\dim \Gamma=d$. Suppose
that $d_i=0$ for some $i =1,\dotsc,k$.  Then for each $\Gamma=
\Gamma_1 \oplus \cdots \oplus \Gamma_k \in \calM \calS(\sigma)$ we
have $\dim \Gamma_i=0$, hence $w(\xi)=0$.  Assume now that $d_i \geq
1$ for $i=1,\dotsc,k$. We have then
\begin{equation}\label{E:weightofxi}
w(\xi)=
\sum_{\Gamma \in  \calM \calS(\sigma) \, : \, \dim \Gamma=d} MV_d(\Gamma_1,\dotsc,\Gamma_k ; (\dim \Gamma_1,\dotsc, \dim \Gamma_k))
\end{equation}
where the sum is taken over all polytopes $\Gamma=\Gamma_1 \oplus \cdots \oplus
\Gamma_k \in \calM \calS(\sigma)$ with $\dim \Gamma=d$ and $\dim \Gamma_i > 0$ for $i=1,\dotsc,k$.
Let $\underline{t}=(t_1,\dotsc,t_k)$ be any collection of positive
integer numbers such that $t_1+\cdots+t_k=d$. Consider the polytope
$$\underbrace{\sigma_1+\cdots+\sigma_1}_{t_1}+\cdots+\underbrace{\sigma_k+\cdots+\sigma_k}_{t_k}$$
together with the convex mixed subdivision induced by the functions
$\tilde{\nu_1}, \dotsc, \tilde{\nu_k}$, where $\tilde{\nu_i}$ is used
for each copy of $\sigma_i$. This mixed subdivision consists of the polytopes
$$\underbrace{\Gamma_1+\cdots+\Gamma_1}_{t_1}+\cdots+\underbrace{\Gamma_k+\cdots+\Gamma_k}_{t_k}$$
for $\Gamma=\Gamma_1 \oplus \cdots \oplus \Gamma_k \in \calM \cal S(\sigma)$.
By Formula~\eqref{E:nonpure} in Lemma~\ref{L:mixedvolumesum}, we get
$$
MV_d(\sigma_1,\dotsc,\sigma_k ; \underline{t})
=
\sum_{\Gamma \in \calM\calS(\sigma)}MV_d(\Gamma_1,\dotsc,\Gamma_k ; \underline{t}).
$$
This sum can actually be taken over all $\Gamma \in
\calM\calS(\sigma)$ such that $\dim \Gamma=d$ and $\dim \Gamma_i \geq
1$ for $i=1,\dotsc,d$ since otherwise $MV_d(\Gamma_1,\dotsc,\Gamma_k ;
\underline{t})=0$.  But now if $\underline{t} \neq (\dim
\Gamma_1,\dotsc,\dim\Gamma_k)$ and $t_1+\cdots+t_k=\dim
\Gamma_1+\cdots+\dim \Gamma_k$, then there exists $i \in
\{1,\dotsc,k\}$ such that $\dim \Gamma_i < t_i$ and thus $\vol_{t_i}
(\Gamma_i)=0$, which implies that $MV_d(\Gamma_1,\dotsc,\Gamma_k ;
\underline{t})=0$. Therefore,

$$
MV_d(\sigma_1,\dotsc,\sigma_k ; \underline{t})
=
\sum_{\Gamma \in \calM\calS(\sigma) \, : \, \dim \Gamma_i=t_i \, \mbox{\tiny for all} \; i} MV_d(\Gamma_1,\dotsc,\Gamma_k ; \underline{t}). $$
and thus
$$
\sum_{ \underline{t} \, : \, t_1+\cdots+t_k=d} MV_d(\sigma_1,\dotsc,\sigma_k ; \underline{t})
=
\sum_{\Gamma \in \calM \calS(\sigma) \, : \, \dim \Gamma=d} MV_d(\Gamma_1,\dotsc,\Gamma_k ; (\dim \Gamma_1,\dotsc, \dim \Gamma_k)).
$$ This proves Formula~\eqref{E:nontransversal}. Finally, If $k=d$,
then there is only one collection $\underline{t}=(t_1,\dotsc,t_d)$ of
positive integer numbers such that $t_1+\cdots+t_d=d$, namely
$\underline{t}=(1,1,\dotsc,1)$. Hence, we get
$w(\xi)=MV_d(\sigma_1,\dotsc,\sigma_k ; (1,1,\dotsc,1))=
MV_d(\sigma_1,\dotsc,\sigma_k)$.  \CQFD

Using our weights as intersection multiplicity numbers,
we get a tropical Bernstein theorem directly from Theorem~\ref{T:weights}.

\begin{cor}
Suppose tropical hypersurfaces $Z_1,\dotsc,Z_n \subset N_{\R} \simeq \R^n$ with Newton polytopes $\Delta_1,\dotsc,\Delta_n$
intersect in finitely many points. Then the total number of intersection points counted with multiplicities is equal to the mixed volume
$MV_n(\Delta_1,\dotsc,\Delta_n)$.
\end{cor}
\begin{proof}
The common intersection points are in one-to-one correspondence with
the polytopes $\sigma=\sigma_1+\cdots+\sigma_n$ in the dual mixed
subdivision $\calM \calS$ of $\Delta_1+\cdots+\Delta_n$.  Each
intersection point is a cell of codimension $n$, hence by
Formula~\eqref{E:nontransversaltop}, Theorem~\ref{T:weights}, the
intersection multiplicity number of the tropical hypersurfaces at this
point is equal to $\mbox{MV}_n(\sigma_1,\dotsc,\sigma_n)$, where
$\sigma= \sigma_1+\cdots+\sigma_n$ is the corresponding polytope in
the mixed subdivision. Hence the total number of intersection points
counted with multiplicities is $\sum_{\sigma \in \calM \calS}
\mbox{MV}_n(\sigma_1,\dotsc,\sigma_n)$.  But this sum is equal to
$\mbox{MV}(\Delta_1,\dotsc,\Delta_n)$ by Formula~\eqref{E:nonpure},
Lemma~\ref{L:mixedvolumesum}.
\end{proof}

\section{Non degenenerate tropical complete intersections}
\label{S:nondeg}

All the definitions in this section build upon the following definition
of a nonsingular tropical hypersurface in the same way as definitions
in Section~\ref{S:toric} built upon that of a nonsingular complex
hypersurface.

\begin{definition}\label{D:tropicalnonsingular}
A tropical hypersurface is nonsingular if its dual polyhedral
subdivision is a primitive (convex) triangulation, that is, a
triangulation whose all simplices are primitive.
\end{definition}
 
This definition is well-established in the case of tropical plane
curves. In the general case, it can be motivated by the fact that
around a vertex corresponding to a primitive $n$-simplex, a tropical
hypersurface coincides with a tropical hypersurface with Newton
polytope this simplex.  But such a simplex is given by a basis of the
ambient lattice $M$, and identifying $M$ with $\Z^n$ via this basis
identifies the simplex with the standard unit simplex in
$\Z^n$. Hence, up to a basis change of the ambient lattice, a non
singular tropical hypersurface coincides around each vertex with a
tropical linear hyperplane.  Nonsingular tropical hypersurfaces with
a given Newton polytope do not always exist.  The simplest example is
given by the non primitive tetrahedron with vertices $(0,0,0)$,
$(1,0,0)$, $(0,1,0)$ and $(1,1,2)$ in $\R^3$ which meets the lattice
$\Z^3$ at its vertices and has thus no primitive triangulation (see
\cite{Be2}).  Recall that a tropical hypersurface lies in $N_{\R}
\simeq {\R}^n$, which is the tropical torus associated with some
lattice $N$.  Hence, at this point, a tropical hypersurface is in fact
a toric tropical hypersurface.  A primitive (convex) triangulation of
a polytope induces a primitive (convex) triangulation of each of its
faces. Recall that the truncation $f^{\Gamma}$ of a tropical
polynomial $f$ to a face $\Gamma$ of its Newton polytope also defines
a tropical hypersurface in the corresponding tropical torus
$N(\Gamma)_{\R}$.  Hence, in contrast to the complex case, if $f$
defines a nonsingular tropical hypersurface in the corresponding
tropical torus, then so do automatically all its truncations.
Comparing with the classical definition~\ref{D:classicalnondegenerate}
of a nondegenerate polynomial, this leads to the following
definition.

\begin{definition}
A tropical polynomial is nondegenerate if all its truncations define
nonsingular tropical hypersurfaces in the corresponding tropical tori,
or equivalently, if its dual polyhedral subdivision is a primitive
triangulation.
\end{definition}

Consider now a collection $(f_1,\dotsc,f_k)$ of tropical polynomials in
$\R[x_1,\dotsc,x_n]$, or more generally in $\R[M]$ with $M \simeq \Z^n$.
Let $\Delta_i$ be the Newton polytope of $f_i$.
Define the associated tropical Cayley polynomial $F \in \R[M \oplus
  {\Z}^k]$
by
\begin{equation}\label{E:Cayleypolynomialbis}
F(x,y)=\sum_{i=1}^k y_if_i(x).
\end{equation}
where the operation are the tropical ones. Its Newton polytope is the
associated Cayley polytope $C(\Delta_1,\dotsc,\Delta_k)$.  We have the
following analogue of the classical
definition~\ref{D:nondegenerateclassicalinters}.

\begin{definition}\label{D:nondegeneratetropicalinters}
The collection $(f_1,\dotsc,f_k)$ of tropical polynomials is {\it nondegenerate}
if the associated Cayley polynomial $F$ is nondegenerate which means that
the dual polyhedral subdivision of $C(\Delta_1,\dotsc,\Delta_k)$
is a primitive triangulation.
\end{definition}

Recall that a collection $(\Gamma_i)_{i \in I}$ of faces of
$\Delta_1,\dotsc,\Delta_k$ is called admissible if $I \subset
\{1,\dotsc,k\}$ and $\Gamma_I=\sum_{i \in I} \Gamma_i$ is face a of
$\Delta_I=\sum_{i \in I} \Delta_i$.  The faces of
$C(\Delta_1,\dotsc,\Delta_k)$ are exactly the Cayley polytopes of the
admissible collections $(\Gamma_i)_{i \in I}$.  Since a primitive
triangulation of a polytope induces primitive triangulations of its
faces, it follows that if $(f_1,\dotsc,f_k)$ is nondegenerate, then
for any admissible collection $(\Gamma_i)_{i \in I}$ of faces of
$\Delta_1,\dotsc,\Delta_k$, the collection of tropical polynomials
$(f_i^{\Gamma_i})_{i \in I}$ is also nondegenerate.  For simplicity
denote by $Z_i$ the hypersurface defined by $f_i$.  If $\Gamma_i$ is a
face of $\Delta_i$, we will denote by $Z_{i,\Gamma_i} $ the tropical
hypersurface in $N(\Gamma_i)_{\R}$, or in $N_{\R}$, defined by the
truncation of $f_i$ to $\Gamma_i$.  The next result is the tropical
analogue of Proposition~\ref{P:Cayleytrick}.

\begin{prop}\label{P:complete}
The collection $(f_1,\dotsc,f_k)$ of tropical polynomials is nondegenerate
if and only if for any admissible collection $(\Gamma_i)_{i
\in I}$ of faces of $\Delta_1,\dotsc,\Delta_k$ the hypersurfaces
$Z_{i,\Gamma_i}$ have only transversal intersections each with
intersection multiplicity number $1$.
\end{prop}

\begin{proof}
If $(f_1,\dotsc,f_k)$ is nondegenerate, then the corresponding convex
polyhedral subdivision of the Cayley polytope
$C(\Delta_1,\dotsc,\Delta_k)$ is a primitive triangulation, and thus
the corresponding convex mixed subdivision $\calM \calS$ of
$\Delta=\Delta_1+\cdots+\Delta_k$ is tight. In particular, the mixed
subdivision $\calM \calS$ is pure which means that the hypersurfaces
$Z_1,\dotsc,Z_k$ have only transversal intersections. These
intersections are in one-to-one correspondence with polytopes
$$\sigma=\sigma_1 \oplus \cdots \oplus \sigma_k \in \calM \calS$$
such that $d_i:= \dim \sigma_i \geq 1$ for $i=1,\dotsc,k$. Letting $d=\dim \sigma=d_1+\cdots+d_k$,
the intersection multiplicity number of $Z_1,\dotsc,Z_k$ along the cell $\xi$ dual to $\sigma$
is
\begin{equation}\label{E:wdef}
w(\xi)=\left(\prod_{i=1}^k \mbox{vol}_{d_i}(\sigma_{i}) \right)
\cdot [M(\sigma):M(\sigma_1)+ \cdots+M(\sigma_k)]
\end{equation}
We are going to show that 
\begin{equation}\label{E:wnondef}
\mbox{vol}_{d+k-1} \left( C(\sigma_1,\dotsc,\sigma_k) \right)=\left(\prod_{i=1}^k \mbox{vol}_{d_i}(\sigma_{i}) \right)
\cdot [M(\sigma):M(\sigma_1)+ \cdots+M(\sigma_k)]
\end{equation}
Since $C(\sigma_1,\dotsc,\sigma_k)$ is a primitive simplex,
this will imply that $w(\xi)=1$.  Each $\sigma_i$ is a simplex as well
as $C(\sigma_1,\dotsc,\sigma_k)$.  Thus $\mbox{vol}_{d+k-1} \left( C(\sigma_1,\dotsc,\sigma_k)
\right)$ equals the
absolute value of a $(d+k-1)$-determinant $D$ whose columns are the
coordinates with respect to a basis of
$M(C(\sigma_1,\dotsc,\sigma_k))=M(\sigma) \times \Z^k$ of vectors spanning $C(\sigma_1,\dotsc,\sigma_k)$.
The corresponding determinant taken with respect to a basis of
$\left(M(\sigma_1)+\cdots+M(\sigma_k)\right) \times \Z^k$ is a
determinant $\tilde{D}$ which factors into a product of $k$
determinants $D_1, \dotsc,D_k$. Each factor $D_i$ has size $d_i$ and
is a determinant whose columns are the coordinates with respect to a
basis of $M(\sigma_i)$ of vectors spanning the simplex
$\sigma_i$. The absolute value of $D_i$ is just
$\mbox{vol}_{d_i}(\sigma_i)$. Formula~\eqref{E:wnondef} follows now from
Remark~\ref{R:indexdet}. The same arguments also work for any admissible
$(\Gamma_i)_{i \in I}$ since if $(f_1,\dotsc,f_k)$ is nondegenerate then $(f_i^{\Gamma_i})_{i \in I}$ is nondegenerate too.
This show one implication of Proposition~\ref{P:complete}, let us show the reverse one.

Clearly, if for any admissible
collection $(\Gamma_i)_{i \in I}$ of faces of
$\Delta_1,\dotsc,\Delta_k$ the hypersurfaces $Z_{i,\Gamma_i}$ have
only transversal intersections, then the mixed subdivision $\calM
\calS$ of $\Delta=\Delta_1+\cdots+\Delta_k$ is pure. Consider a full
dimensional polytope in the polyhedral subdivision of
$C(\Delta_1,\dotsc,\Delta_k)$. It may be written as a Cayley polytope
$C(\sigma_1,\dotsc,\sigma_k)$ for some $\sigma=\sigma_1\oplus \cdots \oplus
\sigma_k \in \calM \calS$ with $\dim \sigma=\dim \Delta$.
Set as above $d_i=\dim \sigma_i$ and
$d:=\dim \sigma=\dim \Delta$.  Set $I:=\{i \in
\{1,\dotsc,k\} \, , \, d_i \neq 0 \}$.  Then $\sigma$ is dual to a
cell $\xi$ of the common intersection of the hypersurfaces $Z_i$ for $i \in I$.
Moreover, the intersection multiplicity number between these hypersurfaces along $\xi$ is
$\mbox{vol}_{d_I+k-1}( C(\sigma_i , i \in I))$, where
$C(\sigma_i , i \in I)$ is the Cayley polytope associated with
$\sigma_i$ for $i \in I$ and $d_I$ is the dimension of $\sum_{i \in I}
\sigma_i$. This Cayley polytope lies on the face $C(\Delta_i , i \in
I)$ of $C(\Delta_1,\dotsc,\Delta_k)$.
One can check that
$$ \mbox{vol}_{d+k-1} \left(
C(\sigma_1,\dotsc,\sigma_k)\right)=\mbox{vol}_{d_I+k-1} \left(
C(\sigma_i , i \in I)\right).$$ 
Thus both members are equal to $1$ and
it follows that $ C(\sigma_1,\dotsc,\sigma_k)$ is a primitive simplex.
\end{proof}

\section{Complex and real tropical varieties}
\label{cxrealtropvar}

Complex and real tropical varieties were introduced by Mikhalkin in~\cite{MCounting}.
Here we follow \cite{Be1} and reproduce the definition and notations for
the reader's convenience. Consider a Puiseux series $g=\sum_{r\in R} b_r t^r \in \KK^*$.
Recall that $\val(g)$ is the smallest exponent appearing in $g$ (the usual valuation of $g$),
and that we defined $v(g)$ to be $-\val(g)$ (see Section~\ref{S:standard}).
Define the argument $\arg(g)$ to be the usual argument of the coefficient $b_{\val(g)}$ of the monomial
with smallest exponent. Consider the map 
\[
\begin{array}{rrrl} 
 \VV : &(\KK ^*)^n &\longrightarrow &\RR^n
                \times  (S^1)^n\\
                & z  &  \longmapsto &(v(z_1), \dots,
                v(z_n),\arg(z_1), \dots, \arg(z_n)).
\end{array}
\]
or alternatively
\[
\begin{array}{rrrl} 
 \VC : &(\KK ^*)^n &\longrightarrow & (\CC^*)^n\\
                & z  &  \longmapsto &(e^{v(z_1)+i \arg(z_1)},
\dots, e^{v(z_n)+i  \arg(z_n)})
\end{array}
\]

We will define a complex tropical variety as the topological closure
of the image of a variety in $(\KK ^*)^n$ under either $\VC$ or $\VV$.
We will call both homeomorphic objects a complex tropical variety and
use one and the other in turns depending on the context.  For the rest
of this section, let $f,f_1,\dotsc,f_k$ be polynomials in
$\KK[z_1,\dotsc,z_n]$.  We denote by $Z_f$
the zero set of $f$ in $(\KK ^*)^n$ and by $Y$ the common zero set of
$f_1,\dotsc,f_k$ in $(\KK ^*)^n$.

\begin{definition}
The complex tropical hypersurface $\CC Z_{f,\VC}^{\trop}$ (resp. $\CC Z^{\trop}_{f,\VV}$)
associated with $f$ is the topological closure of the image under $\VC$ (resp. $\VV$) of the hypersurface $Z_f$.
The complex tropical intersection $\CC Y_{\VC}^{\trop}$ (resp. $\CC Y^{\trop}_{\VV}$) associated with $f_1,\dotsc,f_k$ is the
topological closure of the image under $\VC$ (resp. $\VV$) of $Y$.
\end{definition}

A polynomial $\sum c_wz^w \in \KK[z_1,\dotsc,z_n]$ is a real polynomial if
the coefficients $a_r$ of each series $c_\omega = \sum_{r\in R} a_r t^r $ are real.
Assume from now on that $f_1,\dotsc,f_k$ and $f$ are real polynomials.

\begin{definition}
The real tropical hypersurface associated with $f$ is
the intersection of $\CC Z_{f,\VC}^{\trop}$ with $(\RR^*)^n$, 
or alternatively
the intersection of $ \CC Z^{\trop}_{f,\VV}$ with $\RR^n\times
\{0,\pi\}^n$. 
More generally, the real tropical complete intersection associated with $f_1,\dotsc,f_k$ is
the intersection of $\CC Y_{\VC}^{\trop}$ with $(\RR^*)^n$,
or alternatively
the intersection of $\CC Y_{\VV}^{\trop}$ with $\RR^n\times \{0,\pi\}^n$.
\end{definition}
See~\cite{Be1} for pictures of real tropical curves.  The sign of a
Puiseux series $g= \sum_{r\in R} b_r t^r \in \KK^*$ is defined to be
the sign of the coefficient $b_{\val(g)}$ of the monomial with
smallest exponent.

For any $\epsilon=(\epsilon_1,\dotsc, \epsilon_n) \in \{0,1\}^n$,
denote by $\RR(\epsilon)$ the connected component of $(\RR^*)^n$
(called orthant) which consists of all $(x_1,\dotsc,x_n)$ such that
$(-1)^{\epsilon_i} x_i >0$ for $i=1,\dotsc,n$.  We keep the notation
$(\R_+)^n$ for the positive orthant which corresponds to
$\epsilon=(0,\dotsc,0)$.  Denote by $\RR Z_{f,\VC,\epsilon}^{\trop}$
the intersection of $\RR Z^{\trop}_{f,\VC}$ with $\RR(\epsilon)$. If
$\epsilon \in {\{0,1\}}^n$, let $\tilde{\epsilon}$ be the element of
$\{0,\pi\}^n$ defined by $\tilde{\epsilon}_i=\pi \Leftrightarrow
\epsilon_i=1$, and define $\RR Z^{\trop}_{f,\VV,\epsilon} \subset
\R^n$ to be the image of $\RR Z^{\trop}_{f,\VV}\cap (\RR^n \times
\{\tilde{\epsilon}\})$ under the natural identification of $\RR^n
\times \{\tilde{\epsilon}\}$ with $\RR^n$.  If $Z^{\trop}_{f}$ is
nonsingular one can reconstruct $\RR Z^{\trop}_{f,\VC,\epsilon}$ only
from the data of $ Z^{\trop}_f $ and the collection of signs of the
coefficients of $f$ (see \cite{Mikhmath.AG/0403015} p. 25 and 37,
\cite{Vir01}, and \cite{Mik00} Appendix for the case of amoebas).
Consider the tropical hypersurface $ Z^{\trop}_f \subset \RR^n$, the
induced subdivision $\Xi_f$ of $\RR^n$ and the dual subdivision
$\calS_f$ of its Newton polytope $\Delta_f$. Let $\delta_f$ be the
sign distribution at the vertices of $\calS_f$ such that a vertex
$\omega$ is labelled with the sign of the corresponding coefficient
$c_{\omega}$ in $f(z)=\sum c_wz^w \in \KK[z_1,\dotsc,z_n]$.

\begin{lemma}\label{realhyp}
Assume $Z^{\trop}_{f}$ is nonsingular. Then its positive part $\RR
Z^{\trop}_{f,\VV,(0,\dotsc,0)}$ is the closure of the union of the
$(n-1)$-cells of $Z^{\trop}_{f}$ which are dual to edges with vertices
getting different signs via $\delta_f$. More generally, let $\epsilon
\in \{0,1\}^n$ and define the polynomial $f_{\epsilon}$ by
$f_\epsilon(x_1,\cdots,x_n)=f((-1)^{\epsilon_1}x_1,
\dotsc,(-1)^{\epsilon_n}x_n)$.  Then, $\RR Z^{\trop}_{f,\VV,\epsilon}$
is the closure of the union of the $(n-1)$-cells of
$Z^{\trop}_{f_{\epsilon}}$ which are dual to edges with vertices
getting different signs via $\delta_{f_{\epsilon}}$.
\end{lemma}

It is worth noting that $\RR Z^{\trop}_{f,\VV,\epsilon}$ and $\RR
Z^{\trop}_{f,\VC,\epsilon}$ are homeomorphic for each $\epsilon \in
\{0,1\}^n$. In particular $\RR Z^{\trop}_{f,\VV}$ and $\RR
Z^{\trop}_{f,\VC}$ are homeomorphic.  We use the notations of
Section~\ref{S:tcon}. Let $H \subset \Delta^*$ be the piecewise-linear
hypersurface which is constructed by means of the combinatorial
patchworking out of the data $\calS_f$ and $\delta_f$. As a direct
consequence of Lemma~\ref{realhyp}, we obtain the following result.

\begin{prop}\label{P:tropThyper}
Assume that $Z^{\trop}_{f}$ is nonsingular, or equivalently, that the
subdivision $\calS_f$ is a primitive triangulation. Then there exists
an homeomorphism $h:(\RR^*)^n \to (\Delta \setminus \partial\Delta)
^*$ such that $h(\RR Z^{\trop}_{f,\VV})= H \cap (\Delta \setminus
\partial\Delta)^*$. The same property holds for $\RR
Z^{\trop}_{f,\VC}$.
\end{prop}
Here $(\Delta \setminus \partial\Delta)^*$ is the union of the $2^n$
symmetric copies of the relative interior of $\Delta$ under the
hyperplane reflections. Denote by $\Delta_1,\dotsc,\Delta_k$ the
Newton polytopes of $f_1,\dotsc,f_k$, respectively, and set
$\Delta=\Delta_1+\cdots+\Delta_k$. Each polynomial $f_i$ determines a
convex polyhedral subdivision $\calS_i$ of $\Delta_i$ and a sign
distribution $\delta_i$ at the vertices of $\calS_i$. Consider the
piecewise-linear hypersurface $H_i^{\Delta, *} \subset \Delta^*$
constructed out of these data by means of the combinatorial
patchworking for complete intersections (see Section~\ref{S:tcon}).

\begin{prop}\label{P:tropTcomp}
Assume that $f_1,\dotsc,f_k$ define a nondegenerate tropical complete
intersection, which means that the corresponding convex polyhedral
subdivision of the Cayley polytope $C(\Delta_1,\dotsc,\Delta_k)$ is a
primitive triangulation. Then, there exists an homeomorphism
$h:(\RR^*)^n \to (\Delta \setminus \partial\Delta) ^*$ such that
$h(\RR Z^{\trop}_{f_i,\VV})= H_i^{\Delta, *} \cap (\Delta \setminus
\partial\Delta)^*$ for $i=1,\dotsc,k$. The similar property holds for
the real tropical hypersurfaces $\RR Z^{\trop}_{f_i,\VC}$.
\end{prop}
Therefore, $\RR Y^{\trop}_{\VV}$ (resp., $\RR Y^{\trop}_{\VC}$) is homeomorphic to
the common intersection inside $(\Delta \setminus \partial\Delta)^*$
of the piecewise-linear hypersurfaces $H_i^{\Delta, *}$.

Recall that $\CC Y^{\trop}_{\VC}$ is a subset of the torus $(\CC^*)^n$.
We may assume without loss of generality that the polytope $\Delta$ has non empty interior.
Consider the usual compactification of $(\CC^*)^n$ into the toric variety $X_{\Delta}$ associated with $\Delta$,
and let $\iota: (\CC^*)^n \hookrightarrow X_{\Delta}$ denote the corresponding inclusion.
We define the compactification $\overline{\CC Y}^{\trop}_{\VC}$
to be the closure of $\iota(\CC Y^{\trop}_{\VC})$ in
$X_{\Delta}$. 
Note that the stratification of $X_{\Delta}$ into orbits of the action of
$(\CC^*)^n$ defines a natural stratification of $\overline{\CC Y}^{\trop}_{\VC}$.

We sum up natural maps in  the following commutative diagram.
\[
  \xymatrix{
  \RR Y^{\trop}_{\VV}  \ar[r]^\sim\ar@{^{(}->}[d] & \RR Y^{\trop}_{\VC}
   \ar@{^{(}->}[r]\ar@{^{(}->}[d] & {(\RR^{\relax *})^n}\ar@{^{(}->}[d] \ar@{^{(}->}[r]^{\iota_\RR} & \RR X_{\Delta} \ar@{^{(}->}[d] \\
   {\CC Y^{\trop}_{\VV}} \ar[r]^\sim & {\CC Y^{\trop}_{\VC}} \ar@{^{(}->}[r] & (\CC^{\relax *})^n\ar@{^{(}->}[r]^\iota & X_{\Delta}
  }
\]

Define $\overline{\RR Y}^{\trop}_{\VC}$ to be the
intersection of $\overline{{\CC Y}}^{\trop}_{\VC}$ with the
real part $\RR X_{\Delta}$ of $X_{\Delta}$.
Clearly $\overline{\RR Y}^{\trop}_{\VC}$ is also the closure of
$\iota_\RR(\RR Y^{\trop}_{\VC})$ in $\RR X_{\Delta}$. One can see that the natural
stratification of $\overline{\RR Y}^{\trop}_{\VC}$ induced
by the torus action corresponds to the stratification of the
$T$-complete intersection of Theorem~\ref{St} induced by the face complex of $\Delta$.
Consider for $i=1,\dotsc,k$ the piecewise-linear hypersurface 
$\widetilde{H_i^{\Delta}} \subset \widetilde{\Delta}$ (see Section~\ref{S:tcon}).

\begin{prop} Assume that $f_1,\dotsc,f_k$ define a nondegenerate tropical complete intersection.
Then, there exists an homeomorphism $h:\RR X_{\Delta} \rightarrow \widetilde{\Delta}$
sending $\overline{\RR Y}^{\trop}_{\VC}$ to the common intersection of the
piecewise-linear hypersurfaces  $\widetilde{H_i^{\Delta}}$.
\end{prop}

\section{E-polynomials and mixed signature}
\label{S:Epol}
We recall briefly definitions and some properties of so-called
E-polynomials, see ~\cite{Ba-Bo} and~\cite{Da-Kho}.  Let $X$ be a
quasi-projective algebraic variety over $\C$.  For each pair of
integers $(p,q)$, we set
$$
e^{p,q}(X)=\sum_{k \geq 0} {(-1)}^k h^{p,q}(H_c^k(X)),$$
where $h^{p,q}(H_c^ k(X))$ is the dimension of the $(p,q)$-component of the mixed Hodge
structure of the $k$-th cohomology with compact supports. If $X$ is a nonsingular projective variety
then we have $e^{p,q}(X)={(-1)}^{p+q} h^{p,q}(X)$ (see~\cite{Da-Kho}).

The {\it E-polynomial} of $X$ is the sum

\begin{equation}\label{E:Epoly}
E(X;u,v)=\sum_{p,q} e^{p,q}(X)u^pv^q,
\end{equation}
(see~\cite{Ba-Bo}, and~\cite{Da-Kho} where $E(X;u,\bar{u})$ was introduced).
We have the following properties.

\begin{itemize}
\item
If $X$ is a disjoint union of a finite number of locally closed varieties $X_i$, $i \in I$, then
\begin{equation}\label{E:Eadditive}
E(X;u,v)=\sum_{i \in I} E(X_i;u,v).
\end{equation}
\item
\begin{equation}\label{E:Emultiplicative}
E(X \times Y;u,v)=E(X;u,v) \cdot E(Y;u,v).
\end{equation}
\item
If $\pi: Y \rightarrow X$ is a locally trivial fibration with respect to the Zarisky topology
and $F$ is the fiber over a closed point of $X$, then
\begin{equation}\label{E:Efibration}
E(Y;u,v)=E(X;u,v) \cdot E(F;u,v).
\end{equation}
\end{itemize}

In particular, we get (see~\cite{Da-Kho}) $E({\C}P^1;u,v)=1+uv$, $E({\C};u,v)=uv$, $ E({\C}^*;u,v)=uv-1$
and thus
$$E({\C}^k;u,v)=u^kv^k \; , \quad E({({\C}^*)}^k;u,v)=(uv-1)^k
$$

Let us define
\begin{equation}\label{E:varphi}
\varphi(u):= \frac{E(X;1,u)+E(X;u,1)}{2}
\end{equation}
and
\begin{equation}\label{E:extendedsigma}
\tilde{\sigma}(X):=\varphi(-1).
\end{equation}
We will call $\tilde{\sigma}(X)$ the {\bf mixed signature} of $X$.
This is justified by the following result.

\begin{prop}
If $X$ is a nonsingular projective variety then its mixed signature
and usual signature coincide:

$$\tilde{\sigma}(X)=\sigma(X).$$
\end{prop}

\proof
If $X$ is a nonsingular projective variety then
$$\sigma(X)= \sum_{p+q=0 \bmod{2}}{(-1)}^p h^{p,q}(X),$$
where $h^{p,q}(X)$ is the usual Hodge number of type $(p,q)$ of $X$,
and the result follows from the fact that $e^{p,q}(X)={(-1)}^{p+q} h^{p,q}(X)$ (see~\cite{Da-Kho}).
{\CQFD}

The mixed signature of a complex torus is given by

\begin{equation}\label{E:signaturetore}
\tilde{\sigma}({({\C}^*)}^k)={(-2)}^k
\end{equation}

The additivity of the $E$-polynomial implies that of the mixed signature.

\begin{prop}
If $X$ is a disjoint union of a finite number of locally closed varieties $X_i$, $i \in I$, then
$$\tilde{\sigma}(X)=\sum_{i \in I} \tilde{\sigma}(X_i)$$
\end{prop}

Following ~\cite{Da-Kho} we show how the mixed signature of a toric complete intersection
can be expressed in terms of mixed signatures of toric hypersurfaces.

Consider polynomials $f_1,f_2\dotsc,f_k \in \C[x]$, $x=(x_1,\dotsc,x_n)$, which define a toric complete
intersection

$$Y=\{f_1=f_2=\dotsc=f_k=0\} \subset {(\C^*)}^n.$$
Introduce auxiliary coordinates $y_1,\dotsc,y_k$ and for $I \subset \{1,\dotsc,k\}$ define the toric hypersurface $X_I$ by

$$X_I =\left\{\sum_{i \in I} y_if_i(x)-1=0 \right\} \subset{(\C^*)}^{n+\mid I \mid}.$$

\begin{prop}\label{P:sigmaYsigmaX}
We have

$$
\tilde{\sigma}(Y)= {(-2)}^n+{(-1)}^k \sum_{I \subset \{1,\dotsc,n\}} \tilde{\sigma}(X_I).$$
\end{prop}

\proof

Denote by $X$ the hypersurface in ${(\C^*)}^{n} \times {\C}^k$ with
equation $\sum_{i=1}^k y_if_i(x)-1=0$. The restriction to $X$ of the
projection ${(\C^*)}^{n} \times {\C}^k \rightarrow {(\C^*)}^{n}$ is a
locally trivial fibration over ${(\C^*)}^{n} \setminus Y$ with each
fiber a linear subspace of ${\C}^k$. It follows from the properties of
the E-polynomial that

$$
\begin{array}{lll}
E(X;u,v) & = & E({\C}^{k-1}; u,v) \cdot \left[ E({(\C^*)}^{n}; u,v)-E(Y;u,v) \right] \\
  & = & (uv)^{k-1} \cdot \left[ (uv-1)^n-E(Y) \right].
\end{array}
$$
Passing to the mixed signature yields
$$\tilde{\sigma}(X)={(-1)}^{k-1} \left[{(-2)}^n-\tilde{\sigma}(Y)\right].$$
By additivity, we have
$$\tilde{\sigma}(X)=\sum_{I \subset \{1,\dotsc,n\}} \tilde{\sigma}(X_I)$$
and the result follows.{\CQFD}

Assume now that the polynomials $f_1,\dotsc,f_k$ which define $Y$ are
real polynomials so that $Y$ and the hypersurfaces $X_I$ are in turn real
algebraic varieties. We are interested in the (topological) Euler
characteristic of ${\R}Y$. Like the E-polynomial, the Euler
characteristic is additive and multiplicative. We have $\chi({(\R^*)}^{k})={(-2)}^{k}$
and thus comparing with~\eqref{E:signaturetore} we obtain

\begin{equation}\label{E:eulertorus}
\tilde{\sigma}({(\C^*)}^{k})=\chi({(\R^*)}^{k})
\end{equation}

Recall that a toric variety is a real variety (defined by polynomial
equations with real coefficients) and is the disjoint union of
torus orbits. The additivity of the mixed signature and the
Euler characteristic together with Formula~\eqref{E:eulertorus} imply the following result,
which will be not used after.

\begin{prop}
For any toric variety $X$, we have
$$\tilde{\sigma}(X)=\chi({\R}X)$$
\end{prop}

We obtain the following analogue of Proposition~\ref{P:sigmaYsigmaX} for the Euler characteristic
of the real part in place of the mixed signature.

\begin{prop}\label{P:eulerYeulerX}
We have 
$$
\chi({\R}Y)= {(-2)}^n+{(-1)}^k \sum_{I \subset \{1,\dotsc,n\}} \chi({\R}X_I).$$
\end{prop}

\proof

We adapt the proof of  Proposition~\ref{P:sigmaYsigmaX}.
Let $X$ the hypersurface in ${(\C^*)}^{n} \times {\C}^k$ with equation
$\sum_{i=1}^k y_if_i(x)-1=0$. Using the projection ${(\C^*)}^{n} \times {\C}^k \rightarrow {(\C^*)}^{n}$
together with the additivity and multiplicativity properties of the Euler characteristic yields

$$\chi({\R}X)=\chi({\R}^{k-1}) \cdot \left[ \chi({(\R^*)}^{n})-\chi({\R}Y) \right]$$
and thus
$$\chi({\R}X)={(-1)}^{k-1} \left[{(-2)}^n-\chi({\R}Y)\right].$$
By additivity, we have
$$\chi({\R}X)=\sum_{I \subset \{1,\dotsc,n\}} \chi({\R}X_I)$$
and the result follows.{\CQFD}

\section{Statement of the main result}
\label{S:mainresult}
We use the notations of Section~\ref{cxrealtropvar}.
Consider real polynomials $f_1,\dotsc,f_k \in \KK[z_1,\dotsc ,z_n]$
with Newton polytopes $\Delta_1,\dotsc,\Delta_k$. Denote by $Y^{\trop}$
the corresponding tropical intersection in (the tropical torus) $\R^n$.
Let $\RR Y^{\trop}$ denote the real tropical intersection with respect to either the map
$\VV$ or the map $\VC$: $\RR Y^{\trop}=\RR Y_{\VV}^{\trop}$ or $\RR Y^{\trop}=\RR Y_{\VC}^{\trop}$
($\RR Y_{\VV}^{\trop}$ and $\RR Y_{\VC}^{\trop}$ are homeomorphic).
Recall that $Y^{\trop}$ is a nondegenerate tropical complete intersection
if and only if the corresponding convex polyhedral subdivision of the Cayley polytope $C(\Delta_1,\dotsc,\Delta_k)$
is a primitive triangulation.

\begin{thm}\label{maintorus}
Assume that $Y^{\trop}$ is a nondegenerate tropical complete intersection.
\begin{enumerate}
\item
The Euler characteristic of $\RR Y^{\trop}$ depends only on the polytopes $\Delta_1,\dotsc,\Delta_k$
and is equal to the mixed signature of a complete intersection in the complex torus
$(\C^*)^n$ of algebraic hypersurfaces with Newton polytopes $\Delta_1,\dotsc,\Delta_k$, respectively.
In other words, if $Y^{alg}$ denotes such a complete intersection in $(\C^*)^n$, then its mixed
signature $\tilde{\sigma}(Y_{alg})$ depends only on $\Delta_1,\dotsc,\Delta_k$ and we have
$$ \chi(\RR Y^{\trop}) = \tilde{\sigma}(Y_{alg}).$$

\item The Euler characteristic of $\overline{\RR Y}^{\trop}$
depends only on the polytopes $\Delta_1,\dotsc,\Delta_k$
and is equal to the mixed signature of a generic intersection in the projective toric variety $X(\Delta)$
of algebraic hypersurfaces with Newton polytopes $\Delta_1,\dotsc,\Delta_k$, respectively.
In other words, if $\overline{Y}^{alg}$ denotes such a generic intersection in $X(\Delta)$, then
$\tilde{\sigma}(\overline{Y}_{alg})$ depends only on $\Delta_1,\dotsc,\Delta_k$ and we have
$$ \chi(\overline{\RR Y}^{\trop}) = \tilde{\sigma}(\overline{Y}_{alg}).$$
\end{enumerate}
\end{thm}

In the second part of Theorem~\ref{maintorus}, we invoke the genericity in order to ensure that
the intersection of $\overline{Y}^{alg}$
with any complex torus orbit in $X(\Delta)$
is a complete intersection in that torus orbit (this latter intersection is
defined by polynomials whose Newton polytopes are faces of $\Delta_1,\dotsc, \Delta_k$).
\smallskip

\noindent {\it Proof of Theorem~\ref{maintorus}.}
Part (2) follows from part (1) using the stratification by torus orbits
and the additivity of the Euler characteristic and that of the mixed signature.
Now, by Proposition~\ref{P:sigmaYsigmaX} and Proposition~\ref{P:eulerYeulerX}, to prove part (1) it suffices
to prove the case $k=1$, that is, the toric hypersurface case. This is the content of the rest of the paper
(see Theorem~\ref{T:hypersurface}).
\CQFD

\section{Mixed signature of a complex toric hypersurface}
\label{S:mixed}
Let $f$ be a nondegenerate Laurent polynomial with Newton polytope
$\Delta \subset {\R}^{n}$ and assume that $\Delta$ has non empty
interior. Denote by $Z \subset {(\C^*)}^n$
the nonsingular hypersurface defined by $f$.

Let $C \subset {\R}^{{n}+1}$ be the cone with vertex $0$ over $\Delta
\times \{1\} \subset {\R}^n \times \R$.  The set of faces of $C$ with
the order given by the inclusion and the rank function $\rho$ given by
the dimension form an {\it Eulerian poset} that we denote by $P$
(See \cite{Ba-Bo}, Example 2.3).  Hereafter, we refer to~\cite{Ba-Bo}
for detailed definitions.  Taking the dual cones of elements in $P$, we
get the dual poset $P^*$ which is an Eulerian poset with rank function
$\rho^*(z^*)={n}+1-\rho(z)$ and rank $n+1$.  If $x \in P$ is any face of
$C$, then we denote by $[x,\hat{1}]$ the sub-poset of $P$ formed by
all the faces of $C$ having $x$ as a face. This is an Eulerian poset
with rank function $z \mapsto \rho(z)-\rho(x)$ and rank
$n+1-\rho(x)=\rho(C)-\rho(x)$.  The dual poset $[x,\hat{1}]^*$ is an
eulerian poset of rank $n+1-\rho(x)$ and with rank function
$\rho^*(z^*)=n+1-\rho(z)$.

Let $M$ denote the lattice ${\Z}^{n+1}$ in which the cone $C$ has its
vertices and which contains the vertices of $\Delta \times \{1\}$. If
$m \in C \cap M$, define $x(m) \in P$ to be the minimal face of $C$
containing $m$ and $\mbox{deg} \, m$ to be the last coordinate of $m$. Hence,
$m=(m_0,\mbox{deg} \, m)$ for some $m_0 \in \mbox{deg} \, m \cdot \Delta$.  The
following result gives a closed formula for $E(Z;u,v)$ in terms of
so-called {\it B-polynomials} of the sub-posets of $P^*$, which are
defined by induction on the rank (see \cite{Ba-Bo}, Definition 2.7).

\begin{thm}[~\cite{Ba-Bo}, Theorem 3.24]\label{T:EforZ}
$$
E(Z;u,v)=\frac{(uv-1)^n}{uv}+
\frac{(-1)^{n+1}}{uv} \sum_{m \in C \cap M} (v-u)^{\rho(x(m))}B([x(m),\hat{1}]^*;u,v){\left(\frac{u}{v} \right)}^{\text{deg} \, m}.
$$
\end{thm}

Define two functions (see \cite{Ba-Bo}, Definition 3.5)

$$S(C,t):=(1-t)^{n+1} \sum_{m \in C \cap M} t^{\text{deg} \, m}$$
and
$$T(C,t):=(1-t)^{n+1} \sum_{m \in Int(C)\, \cap M} t^{\text{deg} \, m},$$
where $Int(C)$ is the interior of $C$.
They satisfy the duality relation (\cite{Ba-Bo}, Proposition 3.6)
\begin{equation}\label{E:dualityST}
S(C,t)=t^{n+1}T(C,t^{-1})
\end{equation}
In fact, $S(C,t)$ is a polynomial of degree $n$ (see, for example, \cite{Br} or Lemma~\ref{L:coefehrS} below).
The sum
$$\sum_{m \in Int(C) \cap M} t^{\text{deg} \, m}$$
can be written as
$$\sum_{\lambda=0} ^{+\infty} Ehr_{\Delta}(\lambda) t^{\lambda},$$
where $Ehr_{\Delta}(\lambda)$ is the number of integer points in
$\lambda \cdot \Delta$.  The number $Ehr_{\Delta}(\lambda)$ can be
expressed as a polynomial of degree $n=\dim(\Delta)$ in $\lambda$
called the {\it Ehrhart polynomial} of $\Delta$. Let $a_l^{\Delta}$ be
the coefficient of $\lambda^l$ in this polynomial:
$$Ehr_{\Delta}(\lambda)=\sum_{l=0}^{n} a_l^{\Delta} \lambda^l.$$
Let $\psi_i$ be the coefficient of $t^i$ in $S(C,t)$:
$$S(C,t)=\sum_{i=0}^{\infty} \psi_i t^i.$$

The following lemma can be found in Section~4.1 of \cite{DCZ} (see
also \cite{DHZ} p. 233) or \cite{Be1}.
\begin{lemma}\label{L:coefehrS}
  One has
$$\psi_i=
\sum_{l=0}^n \left(\sum_{p=0}^i {(-1)}^{i-p} \bin{n+1}{i-p}p^l \right) a_l^{\Delta}
$$
and $\psi_i=0$ for $i \geq n+1$.
\end{lemma}

We are now able to state our main formula for the mixed signature of toric hypersurfaces.

\begin{prop}\label{P:signaturehyper}
One has
\begin{equation}\label{E:signature}
\tilde{\sigma}(Z) =  -(-2)^{n} + \sum_{l=0}^{n} a^\Delta_l (\sum_{i=0}^{n} \sum_{p=0}^{i} (-1)^{n+p}
\bin{n+1}{i-p} p^l).
\end{equation}
\end{prop}
\proof

Recall that $\tilde{\sigma}(Z)=\varphi(-1)$, and that $\varphi(u)=[E(Z;1,u)+E(Z;u,1)]/2$.
Writing $I_m$ for $[x(m), \hat{1}]^*$ in Theorem~\ref{T:EforZ} yields

\begin{eqnarray*}
\varphi(u) & = & \frac{(u-1)^n}{u}+ \\
           &   & \frac{(-1)^{n+1}}{2u}
\sum_{m \in C \cap M}(u-1)^{\rho(x(m))}
\left[B(I_m;1,u) u^{-\text{deg} \, m} + {(-1)}^{\rho(x(m))}B(I_m;u,1)u^{\text{deg} \, m} \right]
\end{eqnarray*}

From \cite{Ba-Bo}, Definition 2.7 and Proposition 2.10, we have that $B(I_m;1,u)=1$ if $m \in Int(C)$ and $B(I_m;1,u)=0$ otherwise, and
that $B(I_m;u,1)=(1-u)^{n+1-\rho(x(m))}$. It follows that
\begin{eqnarray*}
\varphi(u) & = & \frac{(u-1)^n}{u}+
        \frac{(u-1)^{n+1}}{2u} \sum_{m \in C \, \cap M} u^{\text{deg} \, m}
+
\frac{(1-u)^{n+1}}{2u} \sum_{m \in Int(C) \,  \, \cap M} u^{-\text{deg} \, m}
\end{eqnarray*}

The second and third terms in this sum are easily shown to be equal to $\frac{(-1)^{n+1}}{2u}S(C,u)$ and
$\frac{(-u)^{n+1}}{2u}T(C,u^{-1})$, respectively.
Now, in view of the duality~\eqref{E:dualityST}, this gives
\begin{eqnarray*}
\varphi(u) & = & \frac{(u-1)^n}{u}+
        \frac{(-1)^{n+1}}{u}S(C,u).
\end{eqnarray*}
The result follows then using Lemma~\ref{L:coefehrS} and putting $u=-1$. \CQFD

\section{Euler Characteristic of a real nonsingular tropical toric hypersurface}
\label{S:Euler}
Let $X$ be a real nonsingular tropical hypersurface with Newton polytope
$\Delta \subset {\R}^n$ where $\Delta$ is assumed to have non empty interior.
Hence, the dual polyhedral subdivision $\calS$ of $\Delta$ is a primitive triangulation.

\begin{lemma}[\cite{It97}]\label{L:nonempty}
Consider a $k$-simplex of $\calS$ which is contained in the interior of $\Delta$.
Its number of non empty symmetric copies
is $2^n-2^{n-k}$.
\end{lemma}

Denote by $nb_k^{\Delta}$ the number of $k$-simplices of $\calS$ which
are contained in the interior of $\Delta$.  We will see that these
numbers are in fact independent of the chosen primitive triangulation
of $\Delta$.

Let $S_2$ be the Stirling number of the second kind defined by
$$
S_2(i,j)=\frac{1}{j!} \sum_{t=0}^j
 {(-1)}^{j-t} \bin{j}{t}t^i.$$

\begin{prop}[see \cite{Be1,D}]
\label{P:numberofsimplices}
We have
$$nb_k^{\Delta}=\sum_{l=k}^{n} k! S_2(l+1,k+1) {(-1)}^{n-l} a_{l}^{\Delta}$$
\end{prop}
\begin{prop}\label{P:Euler}
The Euler characteristic of $X$ verifies
$$
\chi(X) =  (-1)^{n} \sum_{k=1}^{n}
\frac{2^n-2^{n-k}}{k+1} \sum_{l = k}^{n} \sum_{t=0}^{k+1} (-1)^{t+l}
\bin{k+1}{t} t^{l+1} a_{l}^\Delta
$$
\end{prop}
\proof

The $(k-1)$-cells in the cellular decomposition of the T-hypersurface
corresponding to $X$ are given by the non-empty symmetric copies of
the $k$-simplices of $\calS$ contained in the interior of $\Delta$.
Hence, this number of $(k-1)$-cells is equal to
$nb_k^{\Delta}(2^n-2^{n-k})$ by Lemma~\ref{L:nonempty}. This gives

$$
\chi(X) =  \sum_{k=1}^{n}{(-1)}^{k-1} nb_k^{\Delta}(2^n-2^{n-k})
.$$
Using Proposition~\ref{P:numberofsimplices}, we obtain
$$
\chi(X) = \sum_{k=1}^{n}{(-1)}^{k-1}
(2^n-2^{n-k})\sum_{l=k}^{n}{(-1)}^{n-l} a_{l}^{\Delta}
\sum_{t=0}^{k+1}{(-1)}^{k+1-t}\bin{k+1}{t}t^{l+1} ,$$
and the result follows.
\CQFD
 
\section{Main result for a toric hypersurface}
\label{S:mainhyper}
Let $X$ be any real nonsingular tropical hypersurface with Newton
polytope $\Delta \subset {\R}^n$. We may assume without loss of generality that $\dim \Delta =n$.
Let $Z \subset {(\C^*)}^n$ be any nonsingular
hypersurface defined by a polynomial with Newton polytope $\Delta$.

\begin{thm}\label{T:hypersurface}
We have
$$\chi(X)=\tilde{\sigma}(Z).$$
\end{thm}

This section is mainly devoted to the proof of Theorem~\ref{T:hypersurface}.
Before we give two technical results that are repeately used.

\begin{lemma} [See \cite{Be1} appendix, and \cite{LiW} p. 71]\label{vanlint}
Let $l$ and  $i$ be
nonnegative integers. Then,
for $l + 1 \le i$, one has
\begin{equation}\label{vl}
\sum_{q=0}^{i} (-1)^q \bin{i}{q} q^{l} \;= \; \sum_{q=0}^{i} (-1)^q
\bin{i}{q}(i-q)^{l} = 0
\end{equation}
and, as a consequence, for any integer $p$,
\begin{equation}\label{nbsurj}
\sum_{q=0}^{i} (-1)^q
\bin{i}{q}(p-q)^{l} = 0.
\end{equation}
\end{lemma}

\begin{lemma} [See~\cite{Be1} appendix]\label{binrmq}
One has
$\sum_{t=0}^{p} 2^{p-t}
\bin{t}{k} = \sum_{l=k+1}^{p+1} \bin{p+1}{l}$.
\end{lemma}

Let us now begin the proof of Theorem~\ref{T:hypersurface}.
From Proposition~\ref{P:signaturehyper}, we have
\begin{equation}\label{E:signaturehyper}
  \tilde{\sigma}(Z) =  -(-2)^{n} + \sum_{l=0}^{n} a^\Delta_l (\sum_{i=0}^{n} \sum_{p=0}^{i} (-1)^{n+p}
  \bin{n+1}{i-p} p^l).
\end{equation} 
On the other hand, Proposition~\ref{P:Euler} tell us that 
\begin{equation}\label{E:Euler}
\chi(X) =(-1)^{n+1} \sum_{k=1}^{n}
\frac{2^n-2^{n-k}}{k+1} \sum_{l = k+1}^{n+1} \sum_{t=0}^{k+1} (-1)^{t+l}
\bin{k+1}{t} t^{l} a_{l-1}^\Delta
\end{equation}
Note that the sum on $l$ can be taken from $1$ (and in fact from $0$)
according to Lemma~\ref{vanlint}.
Write
$$\tilde{\sigma}(Z)=-(-2)^{n} + \sum_{l=0}^{n} S_{l,n} \cdot a^\Delta_l,$$
with
\begin{equation}\label{E:coefsignature}
S_{l,n}=(-1)^n\sum_{i=0}^{n} \sum_{p=0}^{i} (-1)^{p}
\bin{n+1}{i-p} p^l,
\end{equation}
and
$$
\chi(X)= \sum_{l=1}^{n} C_{l,n} \cdot a^\Delta_l$$
with
\begin{equation}\label{E:coefEuler}
C_{l,n} = (-1)^{n-l} \sum_{k=1}^{n}
\frac{2^n-2^{n-k}}{k+1}  \sum_{t=0}^{k+1} (-1)^t
\bin{k+1}{t} t^{l+1}.
\end{equation}

\begin{lemma}\label{L:recurSC}
We have 
\begin{eqnarray*}
S_{l,n+1}& = & -2 S_{l,n} \quad \mbox{if} \quad l \neq 0\\
C_{l,n+1}& = & -2 C_{l,n},
\end{eqnarray*}
and $S_{0,n} =  (-2)^n$.
\end{lemma}

\proof We have
\begin{eqnarray*}
S_{0,n} & = & (-1)^{n}\sum_{i=0}^{n} \sum_{p=0}^{i} (-1)^{p}
\bin{n+1}{i-p}\\
 & = & (-1)^{n}\sum_{i=0}^{n} \sum_{b=0}^{i} (-1)^{i-b}
\bin{n+1}{b}\\
& = & (-1)^{n} \sum_{b=0}^{n} (-1)^{b}
\bin{n+1}{b} \sum_{t=b}^{n} (-1)^t\\
& = & (-1)^{n} \sum_{b=0,\, b=n \bmod{2}}^{n} \bin{n+1}{b}.
\end{eqnarray*}
The sum and difference
$$
\sum_{b=0,\, b=n \bmod{2}}^{n} \bin{n+1}{b} \; \pm \sum_{b=0, \, b=n+1 \bmod{2}}^{n+1} \bin{n+1}{b}$$
are equal to $2^{n+1}$ and $0$, respectively. This yields $S_{0,n} =  (-2)^n$.
We have
$$
S_{l,n+1}=(-1)^{n}\sum_{i=0}^{n+1} \sum_{p=0}^{i} (-1)^{p+1}
\bin{n+2}{i-p} p^l.
$$
Use that $\bin{n+2}{i-p}= \bin{n+1}{i-p-1} + \bin{n+1}{i-p}$ to obtain
$$
(-1)^{n}S_{l,n+1}=\sum_{i=0}^{n+1} \sum_{p=0}^{i} (-1)^{p+1}
\bin{n+1}{i-p-1} p^l +\sum_{i=0}^{n+1} \sum_{p=0}^{i} (-1)^{p+1}
\bin{n+1}{i-p} p^l  .
$$
By Lemma~\ref{vanlint} $\sum_{p =0}^{n+1} (-1)^{p+1}
\bin{n+1}{n+1-p}(p)^{l} = 0 $ since $l \leq n$. We have $(-1)^{p+1}
\bin{n+1}{i-p-1} p^l=0$ if $p=0$ and $l \neq 0$, and $\bin{n}{i-i-1}=0$.
Hence for $l \neq 0$, we get
\begin{eqnarray*}
(-1)^{n}S_{l,n+1}&=&\sum_{i=1}^{n+1} \sum_{p=0}^{i-1} (-1)^{p+1}
\bin{n+1}{i-p-1} p^l +\sum_{i=0}^{n} \sum_{p=0}^{i} (-1)^{p+1}
\bin{n+1}{i-p} p^l  \\
&=&\sum_{j=0}^{n} \sum_{p=0}^{j} (-1)^{p+1}
\bin{n+1}{j-p} p^l +\sum_{i=0}^{n} \sum_{p=0}^{i} (-1)^{p+1}
\bin{n+1}{i-p} p^l, 
\end{eqnarray*}
with the change of index $j=i-1$. This gives the equality 
$S_{l,n+1} = -2 S_{l,n}$ for $l \neq 0$.

Finally, Let us show that $C_{l,n+1} = -2 C_{l,n}$.
We have
\begin{eqnarray*}
C_{l,n+1}(-1)^{n-l+1} & = & \sum_{k=1}^{n+1}
\frac{2^{n+1}-2^{n+1-k}}{k+1}  \sum_{t=0}^{k+1} (-1)^t
\bin{k+1}{t} t^{l+1}\\
& =& \sum_{k=1}^{n}
\frac{2^{n+1}-2^{n+1-k}}{k+1}  \sum_{t=0}^{k+1} (-1)^t
\bin{k+1}{t} t^{l+1} +
\frac{2^{n+1}-1}{n+2}  \sum_{t=0}^{n+2} (-1)^t
\bin{n+2}{t} t^{l+1}\\
& =& 2 \sum_{k=1}^{n}
\frac{2^{n}-2^{n-k}}{k+1}  \sum_{t=0}^{k+1} (-1)^t
\bin{k+1}{t} t^{l+1} 
\end{eqnarray*}
since $\sum_{t=0}^{n+2} (-1)^t
\bin{n+2}{t} t^{l+1}=0$ by Lemma~\ref{vanlint} \CQFD

\begin{lemma}\label{L:SC}
We have $S_{n,n}= C_{n,n}$
\end{lemma}

\proof
\begin{eqnarray*}
 C_{n,n}& = & \sum_{k=1}^{n}
\frac{2^n-2^{n-k}}{k+1}  \sum_{t=0}^{k+1} (-1)^t
\bin{k+1}{t} t^{n+1}\\
& = &  \sum_{k=0}^{n}
\frac{2^n-2^{n-k}}{k+1}  \sum_{t=1}^{k+1} (-1)^t
\bin{k+1}{t} t^{n+1}.
\end{eqnarray*}
Just notice that the two changes of range do not affect the sum.
Then use that $\frac{1}{k+1}\bin{k+1}{t} t^{n+1}= \bin{k}{t-1} t^{n}$
to get 
\begin{eqnarray*}
 C_{n,n}& = & \sum_{k=0}^{n}
(2^n-2^{n-k})  \sum_{t=1}^{k+1} (-1)^t
\bin{k}{t-1} t^{n}\\
& = & \sum_{k=0}^{n}
2^n  \sum_{t=1}^{k+1} (-1)^t
\bin{k}{t-1} t^{n} -  \sum_{k=0}^{n}
2^{n-k}  \sum_{t=1}^{k+1} (-1)^t 
\bin{k}{t-1} t^{n}\\
& = &2^n \sum_{t=1}^{n+1} (-1)^t t^{n}  \sum_{k=t-1}^{n}
 \bin{k}{t-1} - \sum_{t=1}^{n+1} (-1)^t t^{n}  \sum_{k=t-1}^{n}
2^{n-k}  \bin{k}{t-1}\\
& = &2^n \sum_{t=1}^{n+1} (-1)^t t^{n}
 \bin{n+1}{t} - \sum_{t=1}^{n+1} (-1)^t t^{n}  \sum_{m=t}^{n+1} \bin{n+1}{m}
\end{eqnarray*}
by Lemma~\ref{binrmq} and the fact that $\sum_{k=t-1}^{n}
\bin{k}{t-1}= \bin{n+1}{t}$. 
Then, the first term is $0$ by
Lemma~\ref{vanlint} and we get
\begin{eqnarray*}
- C_{n,n}& = &\sum_{t=1}^{n+1} (-1)^t t^{n}  \sum_{m=t}^{n+1} \bin{n+1}{m}\\
& = &
  \sum_{m=1}^{n+1}\sum_{t=1}^{m}(-1)^t t^{n} 
 \bin{n+1}{m}\\
& = &
\sum_{k=1}^{n+1}\sum_{t=1}^{k}(-1)^{t} t^{n} 
 \bin{n+1}{k-t}
\end{eqnarray*}
with the change of indices $k=(t-m)+n+1$. The sums over $k$ and $t$ can actually be taken starting from $0$.
Moreover the sum over $k$ can be taken until $n$ since for $k=n+1$
we get $\sum_{t=0}^{n+1}(-1)^t t^{n} \bin{n+1}{n+1-t}$ which is zero due to Lemma~\ref{vanlint}.
This gives
$$
- C_{n,n}=\sum_{k=0}^{n}\sum_{t=0}^{k}(-1)^{t} t^{n} 
 \bin{n+1}{k-t}.
$$

On the other hand, we have
$$S_{n,n}=(-1)^n\sum_{i=0}^{n} \sum_{p=0}^{i} (-1)^{p}
\bin{n+1}{i-p} p^{n}.$$
Therefore, we get $C_{n,n}= (-1)^{n+1} S_{n,n}$ which is the desired equality for $n$
odd. Suppose now that $n$ is even. Taking $l=n$ in~\eqref{E:coefsignature},
and noting that the first sum can be taken until $i=n+1$ due to Lemma~\ref{vanlint},
we get
\begin{eqnarray*}
(-1)^n S_{n,n}&=&\sum_{i=0}^{n+1} \sum_{p=0}^{i} (-1)^{p}
\bin{n+1}{i-p} p^{n}\\
&=&\sum_{i=0}^{n+1} \sum_{m=0}^{i} (-1)^{i-m}
\bin{n+1}{m} (i-m)^{n}
\label{A:second}\\
&=&\sum_{k=0}^{n+1} \sum_{m=0}^{n+1-k} (-1)^{n+1-k-m}
\bin{n+1}{m} (n+1-k-m)^{n}\\
&=&\sum_{k=0}^{n+1} \sum_{t=k}^{n+1} (-1)^{t-k}
\bin{n+1}{n+1-t} (t-k)^{n}\\
&=&\sum_{k=0}^{n+1} \sum_{t=k+1}^{n+1} (-1)^{k-t}
\bin{n+1}{t} (k-t)^{n}.
\label{A:last}
\end{eqnarray*}
with the successive changes of indices $m=i-p$, $k=n+1-i$ and
$t=n+1-m$ and using that $n$ is even. 
Suming up the second and last formula yields
$$
2 (-1)^n S_{n,n}=\sum_{k=0}^{n+1} \sum_{t=0}^{n+1} (-1)^{k-t}
\bin{n+1}{t} (k-t)^{n}$$
which is zero by Lemma~\ref{vanlint}. \CQFD 
\smallskip

\noindent {\bf Proof of Theorem~\ref{T:hypersurface}}
We have $S_{0,n}={(-2)}^n$ by Lemma~\ref{L:recurSC} and
$a_0^{\Delta}=1$ by definition of the Ehrhart polynomial.  Hence,
Formula~\eqref{E:coefsignature} can be written as
$\tilde{\sigma}(Z)=\sum_{l=1}^{n} S_{l,n} \cdot a^\Delta_l$.
Comparing with Formula~\eqref{E:coefEuler}, it remains to prove that
$S_{l,n}=C_{l,n}$ for $l=1,\dotsc,n$.  But this clearly follows from
Lemma~\ref{L:recurSC} and Lemma~\ref{L:SC}. \CQFD

\providecommand{\bysame}{\leavevmode\hbox to3em{\hrulefill}\thinspace}
\providecommand{\MR}{\relax\ifhmode\unskip\space\fi MR }
\providecommand{\MRhref}[2]{%
  \href{http://www.ams.org/mathscinet-getitem?mr=#1}{#2}
}
\providecommand{\href}[2]{#2}

\end{document}